\documentclass[12pt]{amsart}
\usepackage{mathrsfs}
\usepackage{amsfonts}
\usepackage{amsthm}
\newtheorem{theorem}{Theorem}
\newtheorem{lemma}{Lemma}

\newtheorem{proposition}{Proposition}

\begin{document}

\title[A sufficient condition for global regularity]{A sufficient
condition for global regularity of the $\overline{\partial}$-Neumann
operator}

\author{Emil J. Straube}
\address{Department of Mathematics\\
Texas A\&M University\\
College Station, Texas, 77843--3368}
\email{straube@math.tamu.edu}

\thanks{Research supported in part by NSF grants DMS 0100517 and DMS 0500842, and by ESI Vienna}

\date{October 14, 2005}

\begin{abstract} 
A theory of global regularity of the $\overline{\partial}$-Neumann
operator is developed which unifies the two principal approaches to
date, namely the one via compactness due to Kohn-Nirenberg
\cite{KohnNirenberg65} and Catlin \cite{Catlin84b} and the one via
plurisubharmonic defining functions and/or vector fields that commute
approximately with $\overline{\partial}$ due to Boas and the author
\cite{BoasStraube91, BoasStraube93}. 

\end{abstract}

\maketitle

\section{Introduction}
The $\overline{\partial}$-Neumann problem and its regularity theory
play important roles both in several complex variables and in partial
differential equations. In several complex variables, the
$\overline{\partial}$-Neumann problem is intimately connected with
solving the $\overline{\partial}$-equation and with the Bergman
projection; in partial differential equations, it provides a
prototype for an elliptic operator with non-coercive boundary
conditions and (in the case of domains of finite type) for a
subelliptic problem. We refer the reader to the surveys
\cite{BoasStraube99, Christ99, DAngeloKohn99, FuStraube99, McNeal05} and the
monographs \cite{FollandKohn72, ChenShaw01, LiebMichel02} for background material. 

Denote by $\Omega$ a smooth bounded pseudoconvex domain in
$\mathbb{C}^{n}$. For $1 \leq q \leq n$, the complex Laplacian
$\Box_{q}$ is given by $\overline{\partial}^{*}\partial +
\partial\overline{\partial}^{*}$ on
$\mathcal{L}^{2}_{(0,q)}(\Omega)$, the usual Hilbert space of
$(0,q)$-forms with coefficients in $\mathcal{L}^{2}(\Omega)$.
$\Box_{q}$ is self-adjoint and onto, hence has a (self-adjoint)
bounded inverse. This inverse is the $\overline{\partial}$-Neumann
operator $N_{q}$. We say that $N_{q}$ is \emph{globally regular} if
it maps $C^{\infty}_{(0,q)}(\overline{\Omega})$, the Fr\'{e}chet space of
$(0,q)$-forms with coefficients in $C^{\infty}(\overline{\Omega})$ (necessarily
continuously) into itself. We say that $N$ is \emph{exactly regular}
when it maps the $\mathcal{L}^{2}$-Sobolev spaces
$W^{s}_{(0,q)}(\Omega)$ of forms with coefficients in $W^{s}(\Omega)$
to themselves (for $s \geq 0$). Standard embedding theorems show that exact regularity
implies global regularity. (It is rather intriguing that so far in all
cases where global regularity is known, it is actually established
via exact regularity.) 

Kohn and Nirenberg proved in \cite{KohnNirenberg65} that for a class
of operators defined by a quadratic form, which includes the
$\overline{\partial}$-Neumann operator, a so called compactness
estimate implies exact regularity, but they did not
address the question of when such an estimate holds. Catlin then
verified in \cite{Catlin84b} that in the case of the
$\overline{\partial}$-Neumann operator, this approach provides indeed
a viable route to global regularity, by showing that a large class of
domains, defined by a geometric condition, satisfies the requisite
estimate. In addition, Catlin's work provides a general sufficient
condition of a potential theoretic nature for compactness. This
condition was systematically investigated by Sibony \cite{Sibony87b}
(see also his survey \cite{Sibony91}). In particular, Sibony's work
gives examples of domains whose boundaries contain large sets (in the
sense of surface measure) of points of infinite type, yet whose
$\overline{\partial}$-Neumann operator nevertheless satisfies a
compactness estimate (and hence is exactly regular). In
\cite{Takegoshi91}, Takegoshi presented an approach that places a certain boundedness condition on the gradients of the functions, rather than on the functions themselves (as had been the case in Catlin's work). In that sense, it may be viewed as a precursor to \cite{McNeal02}, where McNeal introduced a relaxed version of Catlin's condition based on having uniform bounds on the gradients in the metric induced by the complex Hessian of the functions. In \cite{HeferLieb00}, compactness of the $\overline{\partial}$-Neumann problem is studied from the point of view of solution kernels for $\overline{\partial}$, while \cite{Harrington04} contains results in the spirit of Oka's lemma. Recently, the author gave a simple geometric condition, on domains in $\mathbb{C}^{2}$, that implies compactness. Its relation to the potential theoretic conditions discussed here is not understood at present. For a survey on compactness, we refer the reader to \cite{FuStraube99}.

In \cite{BoasStraube91, BoasStraube93}, Boas and the author presented
a new technique for proving Sobolev estimates for the
$\overline{\partial}$-Neumann operator based on the existence of
families of vector fields that have certain \emph{approximate}
commutator properties with $\overline{\partial}$. In particular, such
families of vector fields exist, and hence the
$\overline{\partial}$-Neumann operator is globally regular, when the
domain $\Omega$ admits a defining function that is plurisubharmonic at
boundary points (that is, its complex Hessian is positive
semi-definite at points of the boundary). This covers for example all
smooth convex domains. Other examples of domains where the existence of these families of vector fields has been verified include domains with circular symmetry (\cite{Chen89}), domains whose boundary is of finite type except for a flat piece that is `nicely' foliated by Riemann surfaces (\cite{StraubeSucheston03}, \cite{ForstnericLaurent05} ), and domains whose weakly pseudoconvex directions at boundary points are limits, from inside, of weakly pseudoconvex directions of level sets of the boundary distance (\cite{Straube01}; this class includes domains whose closure admits a particularly nice Stein neighborhood basis). In \cite{BoasStraube93}, the authors studied in
detail the situation when the weakly pseudoconvex boundary points are
contained in a submanifold $M$ of the boundary having the property
that its real tangent space at each point is contained in the complex
tangent space to the boundary at the point (this happens, for
example, for complex submanifolds of the boundary). They identified a
De Rham cohomology class on $M$ as the (only) obstruction to the
existence of the family of vector fields required in their technique.
In particular, when $M$ has trivial first De Rham cohomology (for
example, when $M$ is simply connected), the
$\overline{\partial}$-Neumann operator on $\Omega$ is exactly
regular. In the case of a complex submanifold $M$ of the boundary,
this cohomology class had appeared earlier in
\cite{BedfordFornaess78} in the context of deciding whether or not the
closure
of $\Omega$ admits a Stein neighborhood basis. Its appearance in
connection with global regularity explains why the critical annulus
in the boundary of the worm domains prevents global regularity
\cite{Barrett92, Christ96, Christ99}, while an annulus in the
boundary of certain other Hartogs domains does not do so
\cite{BoasStraube92}, and why an analytic disc is always benign
\cite{BoasStraube92}. 

More recently, Sucheston and the author showed
(\cite{StraubeSucheston02}, Theorem on page 250) that the conditions that appear in \cite{BoasStraube91, BoasStraube93} (i.e. existence of a defining function plurisubharmonic at points of the boundary, existence of a family of vector fields with suitable approximate commutator properties with $\overline{\partial}$, and vanishing of a
cohomology class on certain submanifolds of the boundary) can be
modified in a natural way so as to become equivalent (and
still imply exact regularity). 

The present paper provides a general sufficient condition for exact regularity. It is trivially satisfied for $(0,q)$-forms when there is a compactness estimate (at the level of $(0,q)$-forms).  Modulo classical results, it is also easily seen to be satisfied for all $q \geq 1$ when the assumptions from \cite{BoasStraube91, BoasStraube93}, in the more general form given in \cite{StraubeSucheston02}, hold. In fact, our condition has a potential theoretic flavor, and it will be seen that the approach in \cite{BoasStraube91, BoasStraube93, StraubeSucheston02} arises from extracting the geometric content of the condition. It is noteworthy that the condition discriminates among the form levels, and that it passes from $(0,q)$-forms to $(0,q+1)$-forms (see Lemma \ref{percol} below). When $q>1$, it is satisfied when $\Omega$ admits a defining function whose complex Hessian has the property that the sum of any $q$ eigenvalues is nonnegative. Thus, in the context of pseudoconvex domains, the recent regularity results in \cite{Herbig-McNeal05} are also covered.

The remainder of the paper is organized as follows. In section 2, we
state our new sufficient condition for exact regularity; this is the main result. In section 3, we show that under the assumptions in Theorem \ref{main}, commutators of certain vector fields with $\overline{\partial}$ and with $\overline{\partial}^{*}$, respectively, are benign, in a technical sense needed in the proof of Theorem \ref{main}. This proof is given in section 4. In section 5 we explain why the  assumptions in Theorem \ref{main} are  satisfied under the conditions in \cite{BoasStraube91,BoasStraube93, StraubeSucheston02}, in particular, when there is a family of vector fields that has good approximate commutator properties with $\overline{\partial}$. Section 6 contains some estimates, also required in section 4 in the proof of Theorem \ref{main}, for operators obtained from the $\overline{\partial}$-Neumann problem by elliptic regularization.

\section{A sufficient condition for global regularity}

For summations over multiindices, a superscript $^{\prime}$
indicates that the summation is over increasing tuples only. For $s$
real, $\| u \|_{s}$ denotes the norm in $W_{(0,q)}^{s}(\Omega)$. (When
$s=0$, we will omit the subscript.) 

\begin{theorem}\label{main}
Let $\Omega$ be a smooth bounded pseudoconvex domain in
$\mathbb{C}^{n}$, $\rho$ a defining function for $\Omega$. Let $1 \leq q \leq n$. Assume that
there is a
constant $C$ such that for
all $\epsilon > 0$ there exist a defining function $\rho_{\epsilon}$
for $\Omega$ and a constant $C_{\epsilon}$ with
\begin{equation}\label{i} 
1/C < | \nabla \rho_{\epsilon}| < C \ \ \ \mathrm{on}\ b\Omega\ , 
\end{equation}
and
\begin{equation}\label{ii} 
\left \|\sum_{|K|=q-1}^{\prime} \left (  \sum_{j,k =1}^{n}
\frac{\partial^{2} \rho_{\epsilon}}{\partial z_{j} \partial
\overline{z_{k}}} \frac{\partial \rho}{\partial \overline{z_{j}}}
\overline{u_{kK}} \right )d\overline{z_{K}} \right \|^{2} \leq
\epsilon \left (\|
\overline{\partial} u\|^{2} + \| \overline{\partial}^{*} u \|^{2}
\right ) + C_{\epsilon} \|u \|_{-1}^{2} 
\end{equation}
for all $u \in C_{(0,q)}^{\infty}(\overline{\Omega}) \cap
dom(\overline{\partial}^{*})$. Then the $\overline{\partial}$-Neumann
operator $N_{q}$ on $(0,q)$-forms is exactly regular in Sobolev
norms, that is
\begin{equation}\label{Sobolev}
\| N_{q} u \|_{s} \leq C_{s} \| u \|_{s} \ ,
\end{equation}
for $s \geq 0$ and all $u \in W_{(0,q)}^{s}(\Omega)$.
\end{theorem}

The specific form of the factor $\epsilon$ in the first term on the
right hand side of \eqref{ii} is not relevant; as long as there is a factor $\sigma(\epsilon)$ with $\sigma(\epsilon) \rightarrow 0$ as $\epsilon \rightarrow 0$, one can always suitably rescale the family of defining functions. In particular, the value of the constant in front of $\epsilon$ is immaterial.

The simplest situation in Theorem \ref{main} occurs
when there is one defining function, say $\rho$, that works for
all $\epsilon$. This covers the case when $N_{q}$ is compact, as well
as the situation considered in \cite{BoasStraube91} when $\Omega$
admits a defining function that is plurisubharmonic at boundary
points, or, as in \cite{Herbig-McNeal05} when $q>1$, a defining function whose complex Hessian at boundary points has the property that the sum of any $q$ eigenvalues is nonnegative. We will show in section 5 that the `vector field method' from \cite{BoasStraube93} is also covered and that in fact, more generally, the (equivalent) sufficient conditions given in \cite{StraubeSucheston02} imply the one in Theorem \ref{main}.

When $N_{q}$ is compact, take any defining function $\rho$ and
set $\rho_{\epsilon} = \rho$ for all $\epsilon$. Observe that then
the left hand side of \eqref{ii}
is bounded by $\| u\|^{2}$, independently of $\epsilon$, which in
turn can be bounded by the right hand side if $N_{q}$ is compact (see
for example \cite{FuStraube99}, Lemma 1.1). 

Now assume that there is a defining function $\rho$ whose complex Hessian is positive semidefinite at boundary points. Then \eqref{ii} holds for $q=1$ (hence, in view of Lemma \ref{percol} below, for all $q \geq 1$). Namely, there is a constant $C$ such
that near the boundary, the complex Hessian of $\rho$ is bounded
below by $C\rho |u|^{2}$ or, equivalently, adding the form minus
$C\rho |u|^{2}$ produces a positive semi-definite form. Applying the
Cauchy-Schwarz inequality to this form pointwise and then integrating
shows that the left hand side of \eqref{ii} is dominated by 
\begin{equation}\label{hessian}
\int_{\Omega} \sum_{j,k =1}^{n}
\frac{\partial^{2} \rho}{\partial z_{j} \partial\overline{z_{k}}}
u_{j}\overline{u_{k}} + \|\sqrt{(-\rho)} u\|^{2} + \|u\|_{V}^{2} \ , 
\end{equation}
where $V$ is a relatively compact subdomain of $\Omega$. The complex
Hessian of a defining function acts as a subelliptic
multiplier of order $1/2$ on $1$-forms (\cite{D'Angelo93},
Proposition 4, Section 6.4.2). Applying this to the sum in the
first term of \eqref{hessian} and combining the result with
interpolation of Sobolev norms shows that this term is dominated by 
$\epsilon(\|\overline{\partial} u \|^{2} +
\|\overline{\partial}^{*} u \|^{2}) + C_{\epsilon}\| u \|_{-1}^{2}$.
$\|\sqrt{(-\rho)} u\|^{2}$ is dominated by $\epsilon \|u\|^{2} +
\|u\|^{2}_{\Omega_{\epsilon}} \leq
\epsilon(\|\overline{\partial}u\|^{2} +
\|\overline{\partial}^{*}u\|^{2}) + \|u\|^{2}_{\Omega_{\epsilon}}$,
where $\Omega_{\epsilon}$ denotes the points $z \in \Omega$ with
$\rho(z) < -\epsilon$. Because of interior elliptic regularity of 
$\overline{\partial} \oplus \overline{\partial}^{*}$,
$\|u\|^{2}_{\Omega_{\epsilon}} + \|u\|_{V}^{2}$ can be estimated from
above by $\epsilon(\|\overline{\partial}u\|^{2} +
\|\overline{\partial}^{*}u\|^{2}) + C_{\epsilon}\|u\|_{-1}^{2}$. Thus
\eqref{hessian} is indeed dominated by the right hand side of
\eqref{ii}.

When $q>1$, the following equivalent reformulation of condition \eqref{ii} is useful. Define the quadratic form $H_{\rho,q}(u,\overline{u})$ by
\begin{equation}\label{h-q}
 H_{\rho,q}(u,\overline{u}) = \sum_{|K|=q-1}^{\prime} \left ( \sum_{j,k =1}^{n}
\frac{\partial^{2} \rho}{\partial z_{j} \partial
\overline{z_{k}}} u_{jK}\overline{u_{kK}} \right )\; .
\end{equation}
We have
\begin{lemma}\label{equiv}
Let $\Omega$ be a smooth bounded pseudoconvex domain, $\rho$ a defining function for $\Omega$, let $1\leq q \leq n$, and let $C$ be a constant. Then, modulo rescaling, a family of defining functions $\rho_{\epsilon}$ satisfies \eqref{i} and \eqref{ii} if and only if it satisfies \eqref{i} and
\begin{equation}\label{iii}
\sup_{\beta \in C^{\infty}_{(0,q-1)}(\overline{\Omega}),
\|\beta\| \leq 1}\left \{ \left |\int_{\Omega}H_{\rho_{\epsilon},q}(\overline{\partial}\rho \wedge \beta,\overline{u}) \right |^{2}\right \} \leq \epsilon \left (\|\overline{\partial}u\|^{2} + \|\overline{\partial}^{*}u\|^{2}
\right ) + \widetilde{C_{\epsilon}}\|u\|_{-1}^{2} \; .
\end{equation}
\end{lemma}

The lemma is obvious when $q=1$. When $q>1$, note that both in $H_{\rho,q}$ and in inner products between $(q-1)$-forms, we may sum over all multi-indices $K$ of length $(q-1)$ and then divide by $(q-1)!$. The (pointwise) inner product of the form on the left hand side of \eqref{ii} with a $(q-1)$-form $\beta = \sum_{|K|=q-1}^{\prime} b_{K}d\overline{z_{K}}$ equals $H_{\rho,q}(\overline{\partial}\rho\wedge \widetilde{\beta}, \overline{u})$, where $\widetilde{\beta} = \sum_{|K|=q-1}^{\prime} \overline{b_{K}}d\overline{z_{K}}$, modulo terms containing a factor $(\overline{\partial\rho /z_{k_{s}})u_{kk_{1}\cdots k_{s}\cdots k_{q-1}}}$ for some $s, 1\leq s \leq q-1$ (replacing $(\partial\rho/\partial\overline{z_{j}})\overline{\beta_{K}}$ by $(\overline{\partial}\rho\wedge\widetilde{\beta})_{jK}$ makes an error of the indicated form). Upon summation over $k_{s}$, these terms give rise to coefficients of the normal part of $u$. Thus their Sobolev-1 norm is bounded by $\|\overline{\partial}u\| + \|\overline{\partial}^{*}u\|$ (see for example the argument in the proof of Lemma \ref{dbar} below, in particular \eqref{dbar6}), and interpolation between Sobolev norms shows that the contribution to the integral over $\Omega$ coming from these terms can be bounded by the right hand side of \eqref{ii} or \eqref{iii}, respectively. Integrating over $\Omega$ and taking the supremum over $\beta$ with $\|\beta\| \leq 1$ thus proves the lemma.
$\Box$

Assume now that there is a defining function $\rho$ with the property that at boundary points, the sum of any $q$ eigenvalues of its complex Hessian is nonnegative. This implies that $H_{\rho,q}$ is positive semidefinite at points of the boundary (in fact, the two properties are equivalent). We take $\rho_{\epsilon}=\rho$. Now the argument proceeds as in the case of a plurisubharmonic defining function, but with the complex Hessian $H_{\rho,1}$ replaced by $H_{\rho,q}$. (Proposition 4, Section 6.4.2 in \cite{D'Angelo93} is only formulated for $H_{\rho,1}(u,\overline{u})$, i.e. for $(0,1)$-forms, but it is easily seen to be true for $H_{\rho,q}(u,\overline{u})$ on $(0,q)$-forms, by considering the auxiliary forms $v_{K}:=\sum_{j}u_{j,K}d\overline{z_{j}}\,$; compare also the proof of Lemma \ref{percol} below.) This establishes \eqref{iii} (with $\rho_{\epsilon}=\rho$). Since \eqref{i} is trivially satisfied, Lemma \ref{equiv} shows that the assumptions in Theorem \ref{main} do indeed hold. 

Condition \eqref{ii} is of a \emph{potential theoretic} flavor. This is not surprising: global regularity of the $\overline{\partial}$-Neumann operator probably cannot be characterized in terms of purely geometric conditions on the boundary (in contrast to the much stronger property of subellipticity, which is characterized by the geometric notion of finite type). Nonetheless, it is interesting to see how to extract a geometric condition from condition \eqref{ii}, and what the result of doing so is. Let
$q=1$ for simplicity. Because $u \in dom(\overline{\partial}^{*})$,
the vector $(u_{1}, \cdots , u_{n})$ formed from the components of
$u$ is complex tangential at the boundary. The left hand side of
\eqref{ii} is thus the square of the $\mathcal{L}^{2}$-norm of the
following quantity: the mixed (complex tangential unit vector -
complex normal unit vector) term in the complex Hessian of
$\rho_{\epsilon}$ times $|u|$. The square of this
$\mathcal{L}^{2}$-norm should be bounded by the right hand side.
Since $\|u\|^{2} \leq C(\|\overline{\partial}u\|^{2} +
\|\overline{\partial}^{*}u\|^{2})$, this suggests that one require
that the mixed term in the Hessian of $\rho_{\epsilon}$ should be a
multiplier in $\mathcal{L}^{2}(\Omega)$ with operator norm of order
$\epsilon$. However, this operator norm is given by the sup-norm of
the multiplier, so that the requirement becomes that this mixed term
be uniformly small of order $\epsilon$. Actually, it suffices that
this be the case at points of the boundary (since compactly supported
terms are under control) and in weakly pseudoconvex directions (since
components of $u$ in strictly pseudoconvex directions are under
control, see section 5 below). This geometrization scheme therefore leads 
precisely to  \cite{BoasStraube91, BoasStraube93, StraubeSucheston02}.

It is interesting to note that if the assumptions in Theorem
\ref{main} are satisfied at level $q$, then they are satisfied at
levels $q+1, q+2$, etc. Whether exact regularity of the
$\overline{\partial}$-Neumann operator similarly passes from
$q$-forms to $(q+1)$-forms seems to be open. (It is easy to see that compactness and subellipticity do, \cite{McNeal04}.)

\begin{lemma}\label{percol}
Suppose the assumptions in Theorem \ref{main} are satisfied at some
level $q$, where $1 \leq q \leq n-1$. Then they are satisfied at
level $q+1$.
\end{lemma}

\emph{Proof}: Let $u \in C_{(0,q+1)}^{\infty}(\overline{\Omega}) \cap
dom(\overline{\partial}^{*})$,
$u=\sum_{|J|=q+1}^{\prime}u_{J}d\overline{z_{J}}$. For $k=1, \cdots ,
n$, we define $q$-forms $v_{k}$ by
$v_{k}:=\sum_{|K|=q}^{\prime}u_{kK}d\overline{z_{K}}$. Then $v_{k} \in
dom(\overline{\partial}^{*})$: if $|L|=q-1$, then
$\sum_{j=1}^{n}(v_{k})_{jL} (\partial \rho/\partial z_{j}) = \sum_{j=1}^{n} u_{kjL} (\partial \rho/\partial z_{j}) = -\sum_{j=1}^{n} u_{jkL} (\partial \rho/\partial z_{j}) = 0$ on $b\Omega$,
because $u \in dom(\overline{\partial}^{*})$. Computing
$\overline{\partial}^{*}v_{k}$ gives $\overline{\partial}^{*}v_{k}  = -\sum_{|L|=q-1}^{\prime} \sum_{j=1}^{n} (\partial
(v_{k})_{jL}/\partial z_{j})d\overline{z_{L}} = -
\sum_{|L|=q-1}^{\prime} \sum_{j=1}^{n} (\partial
u_{kjL}/\partial z_{j})d\overline{z_{L}} = \sum_{|L|=q-1}^{\prime}
\sum_{j=1}^{n} (\partial u_{jkL}/\partial
z_{j})d\overline{z_{L}}$. Note that the coefficient of
$d\overline{z_{L}}$ is, up to sign, a coefficient of
$\overline{\partial}^{*}u$,
namely that of $d\overline{z_{kL}}$. In particular, the
$\mathcal{L}^{2}$-norm of $\overline{\partial}^{*}v_{k}$ is dominated
by the $\mathcal{L}^{2}$-norm of $\overline{\partial}^{*}u$. Also,
the $\mathcal{L}^{2}$-norm of $\overline{\partial}u$ is dominated by
$\|\overline{\partial}u\|+ \|\overline{\partial}^{*}u\|$. This is
because the components of $\overline{\partial}u$ are expressed in
terms of bar derivatives of components of $u$, and these are
dominated by $\|\overline{\partial}u\|+ \|\overline{\partial}^{*}u\|$
(see e.g. \cite{ChenShaw01}, Section 4.3). Let now
$\{\rho_{\epsilon}\}$ be the family of defining functions that exist
according to the assumptions in Theorem \ref{main} for $q$-forms. The
same family, up to rescaling of $\epsilon$, works for $(q+1)$-forms.
We have
\[
\sum_{|K|=q}^{\prime}\int_{\Omega} \left |  \sum_{j,k =1}^{n}
\frac{\partial^{2} \rho_{\epsilon}}{\partial z_{j} \partial
\overline{z_{k}}} \frac{\partial \rho}{\partial \overline{z_{j}}}
\overline{u_{kK}} \right |^{2}  =  \sum_{(m,\widehat{K})}^{\prime}
\int_{\Omega} \left |\sum_{j,k =1}^{n}
\frac{\partial^{2} \rho_{\epsilon}}{\partial z_{j} \partial
\overline{z_{k}}} \frac{\partial \rho}{\partial \overline{z_{j}}}
\overline{u_{km\widehat{K}}} \right |^{2} \ . 
\]
The summation on the right hand side is over all $m$ and
$(q-1)$-tuples $\widehat{K}$ so that $(m,\widehat{K})$ is an
increasing $q$-tuple. Summing over \emph{all} $m$ and over \emph{all} increasing
$(q-1)$-tuples $\widehat{K}$ (thus increasing the sum), and replacing
$u_{km\widehat{K}}$ by $-u_{mk\widehat{K}} = -(v_{m})_{k\widehat{K}}$,
we see that the right hand side is bounded by
\[ 
\sum_{m=1}^{n} \sum_{|\widehat{K}|=q-1}^{\prime} \int_{\Omega}\left | 
\sum_{j,k} \frac{\partial^{2} \rho_{\epsilon}}{\partial z_{j}
\partial \overline{z_{k}}} \frac{\partial \rho}{\partial
\overline{z_{j}}}(v_{m})_{k\widehat{K}} \right|^{2}
\leq \epsilon \sum_{m=1}^{n} \left (\|
\overline{\partial}v_{m} \|^{2} + \|\overline{\partial}^{*}v_{m}
\|^{2} \right ) + C_{\epsilon} \sum_{m=1}^{n} \|v_{m}\|_{-1}^{2} \ .
\]
In the last estimate, we have used that $v_{m}$ is a $q$-form, and
that the family $\{\rho_{\epsilon}\}$ satisfies \eqref{ii} for
$q$-forms. By what was said above, the right hand side is dominated
by $\epsilon \left (\| \overline{\partial}u \|^{2} +
\|\overline{\partial}^{*}u \|^{2} \right ) + C_{\epsilon}
\|u\|_{-1}^{2}$. This completes the proof of Lemma \ref{percol}. $\Box$

The above argument has benefitted from correspondence with Jeff McNeal concerning the remark preceeding the statement of Lemma \ref{percol}.
\smallskip

\underline{\emph{Remark 1}:} Let again $q=1$ for simplicity. In the discussion in the previous paragraph, it suffices that the mixed term in the complex Hessian of a defining function be what one might call a `compactness multiplier'. That is, if $Y = \sum_{k}Y_{k}(\partial/\partial z_{k})$ denotes a complex tangential field of type $(1,0)$, consider the operator $A_{\rho, Y}$ from $dom(\overline{\partial}) \cap dom(\overline{\partial}^{*})$, provided with the graph norm, to $\mathcal{L}^{2}(\Omega)$ defined by
\begin{equation}\label{Aop}
A_{\rho,Y}(u) := \left (\sum_{j,k =1}^{n}
\frac{\partial^{2}\rho}{\partial z_{j} \partial
\overline{z_{k}}}\frac{\partial \rho}{\partial \overline{z_{j}}}
\overline{Y_{k}}\right )|u| \ , \, \  u \in dom(\overline{\partial}) \cap
dom(\overline{\partial}^{*}) \subseteq
\mathcal{L}^{2}_{(0,1)}(\Omega) \ . 
\end{equation}
If $A_{\rho,Y}$ is compact for all $Y$,then \eqref{ii} holds with $\rho_{\epsilon} = \rho$, for all $\epsilon$. This follows from a lemma
in functional analysis characterizing compact operators, and the fact
that $\mathcal{L}^{2}_{(0,1)}(\Omega)$ embeds compactly into
$W^{-1}_{(0,1)}(\Omega)$, see e.g. \cite{KohnNirenberg65}, Lemma 1.1,
\cite{McNeal02}, Lemma 2.1 (also note that replacing the gradient of a normalized defining function in \eqref{ii} by that of another
defining function does not affect compactness). This suggests that one study functions that produce a compact operator as in \eqref{Aop}. Alternatively, consider the operator
\begin{equation}\label{Bop}
B_{\rho}(u) := \sum_{j,k =1}^{n}
\frac{\partial^{2}\rho}{\partial z_{j} \partial
\overline{z_{k}}}\frac{\partial \rho}{\partial \overline{z_{j}}}
\overline{u_{k}} \ , \, \  u \in dom(\overline{\partial}) \cap
dom(\overline{\partial}^{*}) \subseteq
\mathcal{L}^{2}_{(0,1)}(\Omega) \ .
\end{equation}
Then \eqref{ii} holds (for $q=1$) with $\rho_{\epsilon} = \rho$
for all $\epsilon$ for some defining function $\rho$ if and only
if the operator $B_{\rho}$ is compact (by the same characterization of compact operators quoted above). The form of $B_{\rho}$ suggests that one study sesquilinear forms that produce a compact operator as in \eqref{Bop}. These observations hint at a theory of `compactness multipliers', yet to be developed, modeled after Kohn's theory of subelliptic multipliers (see for example \cite{D'Angelo93}, section 6.4). 

We also note that the discussion concerning the operator $B_{\rho}$ provides a
compactness property considerably weaker than compactness of the
$\overline{\partial}$-Neumann operator which still implies global
regularity. To see that existence of a defining function $\rho$ such that the associated operator $B_{\rho}$ is compact is a considerably weaker property than
compactness of the $\overline{\partial}$-Neumann operator, consider
smooth bounded convex domains. They always admit a defining function
$\rho$ which is plurisubharmonic at points of the boundary (see
\cite{BoasStraube91}), so that the associated operator
$B_{\rho}$ is compact, by what was said above. However, the
$\overline{\partial}$-Neumann operator on $(0,1)$-forms is compact
(if and) only if the boundary contains no analytic disc, see \cite{FuStraube98}. 
\smallskip

\underline{\emph{Remark 2}:} Whether or not a family of defining functions $\{\rho_{\epsilon}\}$
with gradients that are uniformly bounded on $b\Omega$ satisfies
\eqref{ii} is entirely determined by (the interplay of) these
gradients (with $b\Omega$). More precisely: if the family
$\{\rho_{\epsilon}\}$ satisfies the assumptions of Theorem
\ref{main}, and $\{\widetilde{\rho_{\epsilon}}\}$ is a family of
defining functions such that $\nabla\widetilde{\rho_{\epsilon}}(z)=\nabla\rho_{\epsilon}(z)$ for all $z \in b\Omega$ and all
$\epsilon$, then the family $\{\widetilde{\rho_{\epsilon}}\}$ also
satisfies the assumptions of Theorem \ref{main}, possibly after
rescaling. To see
this, write (near $b\Omega$) $\widetilde{\rho_{\epsilon}}= g_{\epsilon}\rho_{\epsilon}$ with $g_{\epsilon}(z)=1$ when $z \in
b\Omega$. Then 
\begin{equation}\label{hessian2}
\frac{\partial^{2}(g_{\epsilon}\rho_{\epsilon})}{\partial
z_{j} \partial \overline{z_{k}}} = \frac{\partial^{2}
g_{\epsilon}}{\partial z_{j} \partial\overline{z_{k}}}\rho_{\epsilon}
+ \frac{\partial g_{\epsilon}}{\partial
z_{j}}\frac{\partial\rho_{\epsilon}}{\partial\overline{z_{k}}} +
\frac{\partial
g_{\epsilon}}{\partial\overline{z_{k}}}\frac{\partial\rho_{\epsilon}}
{\partial z_{j}} + g_{\epsilon}\frac{\partial^{2}
\rho_{\epsilon}}{\partial z_{j}\partial\overline{z_{k}}}\ . 
\end{equation}
$\sum_{k=1}^{n}(\partial\rho_{\epsilon}/
\partial\overline{z_{k}}(z)\overline{u_{k,K}(z)}$ is (the conjugate of) a coefficient of the normal component of $u$, and so has its Sobolev-$1$ norm dominated by $C_{\epsilon}(\|\overline{\partial}u\| + \|\overline{\partial}^{*}u\|)$. Similarly, since $g_{\epsilon} \equiv 1$ on $b\Omega$, the tangential
derivative $\sum_{k=1}^{n}u_{k,K}(z)(\partial
g_{\epsilon}/\partial\overline{z_{k}})(z)$ equals zero on the boundary (the derivative is tangential on the boundary because $u \in dom(\overline{\partial}^{*})$). Consequently, its $1$-norm is also bounded by
$C_{\epsilon}(\|\overline{\partial}u\| + \|\overline{\partial}^{*}u\|)$ (by the same argument as for the normal component). Thus the contributions coming
from the first three terms on the right hand side of \eqref{hessian2}
to the left hand side of \eqref{ii} can all be bounded in the manner required by the right hand
side of \eqref{ii}, essentially by the argument used above (see in particular the proof of Lemma \ref{equiv}). In the contribution coming from the last term in
\eqref{hessian2}, $g_{\epsilon}$ acts as a bounded multiplier near
the boundary, say where $|g_{\epsilon}| < 2$. So modulo compactly
supported terms its contribution to the left hand side of \eqref{ii}
can be bounded by that of the Hessian of $\rho_{\epsilon}$, hence by
the right hand side of \eqref{ii}. The phenomenon discussed in this
paragraph has analogues in \cite{BoasStraube91, BoasStraube93}, where 
plurisubharmonicity of a defining function and good approximate
commutator properties of
vector fields with $\overline{\partial}$, respectively, are only needed at boundary points.
\smallskip

\underline{\emph{Remark 3}:} In \cite{Kohn99}, Kohn gave a qualitative version of the result in \cite{BoasStraube91} in the sense that the level in the Sobolev scale up to which estimates hold is tied to the Diederich-Forn\ae ss exponent (\cite{DiederichFornaess77}). The above discussion of the situation when there is a plurisubharmonic defining function suggests the possibility of such an analysis in our case also. A trivial observation is that to get estimates at a fixed level $k$ in the Sobolev scale, one only needs \eqref{ii} for some $\epsilon = \epsilon(k) > 0$ (see section 4 below).

\smallskip
\section{Inner products involving commutators with $\overline{\partial}$ and with $\overline{\partial}^{*}$}

The proof of Theorem \ref{main} requires estimates on inner products involving commutators, with
$\overline{\partial}$ and with $\overline{\partial}^{*}$, of vector
fields formed from the family of defining functions
$\{\rho_{\epsilon}\}$ given in Theorem \ref{main}. These estimates are given in Lemma \ref{dbar} and Lemma \ref{dbar*} below. We start out with a lemma which makes precise the statement
that bar derivatives and complex tangential derivatives are `benign'
for the $\overline{\partial}$-Neumann problem. We state it in the
form given in \cite{BoasStraube91}.
\begin{lemma}\label{benign}
Let $k \in \mathbb{N}$, and let $Y$ be a
vector field of type $(1,0)$, smooth on $\overline{\Omega}$, with
$Y\rho = 0$ on $b\Omega$. Then there is a constant $C$ such that for
$u = \sum_{|J|=q}^{\prime} u_{J}d\overline{z_{J}} \in
C_{(0,q)}^{\infty}(\overline{\Omega}) \cap
dom(\overline{\partial}^{*})$, we have
\begin{equation}\label{benighn1}
\sum_{j,J} \left \|\frac{\partial u_{J}}{\partial \overline{z_{j}}}
\right \|_{k-1}^{2} \leq C\left ( \|\overline{\partial}u\|_{k-1}^{2}
+ \|\overline{\partial}^{*}u\|_{k-1}^{2} + \|u\|_{k-1}^{2} \right ) \
,
\end{equation}
and 
\begin{equation}\label{benign2}
\left \|Yu\right \|_{k-1}^{2} \leq  C\left (
\|\overline{\partial}u\|_{k-1}^{2} +
\|\overline{\partial}^{*}u\|_{k-1}^{2} + \|u\|_{k-1}\|u\|_{k} \right
) \ .
\end{equation}
\end{lemma}

\emph{Proof}: The proof may be found in \cite{BoasStraube91}, page 83,
formulas (2) and (3), or in \cite{ChenShaw01}, Section 6.2. $\Box$ 
\smallskip

Define $h_{\epsilon}$ (near the
boundary) by $e^{h_{\epsilon}}\rho = \rho_{\epsilon}$, where $\rho$ is
some fixed defining function for $\Omega$ which near $b\Omega$ agrees
with the signed boundary distance. Note that then $|\nabla\rho|
\equiv 1$ near $b\Omega$. Because of \eqref{i}, the functions $h_{\epsilon}$ are bounded on $b\Omega$ independently of $\epsilon$. Therefore, we can choose a family $\{h_{\epsilon}\}_{\epsilon > 0} \in C^{\infty}(\overline{\Omega})$ that is bounded on $\overline{\Omega}$ independently of $\epsilon$ and so that $\rho_{\epsilon} = e^{h_{\epsilon}}\rho$ in $\overline{\Omega} \cap V_{\epsilon}$, where $V_{\epsilon}$ is a neighborhood of $b\Omega$ that depends on $\epsilon$. We denote the inner product in
$\mathcal{L}^{2}_{(0,q)}(\Omega)$ by $(\cdot, \cdot)_{(0,q)}$. Later, it will sometimes be convenient to have differential operators act coefficientwise in special boundary charts (see \cite{FollandKohn72}, page 33, or \cite{ChenShaw01}, pages 129-130). We fix a cover of (a neighborhood of) $b\Omega$ by special boundary charts and an associated partition of unity. However, when working in special boundary charts, we will suppress the cutoff functions and the summation over the charts so as not to additionally complicate the notation.

\begin{lemma}\label{dbar}
Let $1 \leq q \leq n$, and assume $\{\rho_{\epsilon}\}$ is a family of
defining functions as in Theorem \ref{main}, and let
$\{X_{\epsilon}\}$ be a family of smooth vector fields of type
$(1,0)$ so that near $b\Omega$ (possibly depending on $\epsilon$),
$X_{\epsilon}$ agrees with
$e^{-h_{\epsilon}}\sum_{j=1}^{n}(\partial \rho/\partial
\overline{z_{j}}) (\partial/\partial z_{j})$. Let $X_{\epsilon}$ act componentwise, either in Euclidean coordinates or in special boundary charts. Then there are a
constant $A$ and constants $C_{\epsilon, g}$, $0<\epsilon<1$, $g \in C^{\infty}(\overline{\Omega})$, such that
for all $u \in
C_{(0,q)}^{\infty}(\overline{\Omega}) \cap
dom(\overline{\partial}^{*}), v \in
C_{(0,q+1)}^{\infty}(\overline{\Omega}) \cap
dom(\overline{\partial}^{*})$ 
\begin{multline}\label{dbar0}
\left | \left (\left [\  \overline{\partial}, X_{\epsilon} \right ]u,
gv \right )_{(0,q+1)} \right | \leq A\sqrt{\epsilon}\|g\|_{\infty}^{2} \left (
\|\overline{\partial}v \|^{2} + \| \overline{\partial}^{*}v \|^{2} \right ) +
A\sqrt{\epsilon}\| u \|_{1}^{2}  \\
+ C_{\epsilon, g} \left (
\|\overline{\partial}u\|^{2} + \|\overline{\partial}^{*}u\|^{2} +
\|v\|_{-1}^{2} \right ) \ .
\end{multline}

\end{lemma}

Here, $\|g\|_{\infty}$ denotes the $\mathcal{L}^{\infty}$-norm on $\Omega$. $A$ denotes a constant independent of $\epsilon, g$, whereas $C_{\epsilon, g}$ is allowed to depend on both $\epsilon$ and $g$.

\emph{Proof}: We give the proof when $X_{\epsilon}$ acts componentwise in Euclidean coordinates. When it acts in special boundary charts, the change in the commutator with $\overline{\partial}$ contains only terms of order zero and terms involving bar derivatives of $u$ (letting $X_{\epsilon}$ act coefficientwise in special boundary charts changes the operator by a $0$-th order, albeit nonscalar, operator). The contribution from these terms can be estimated by the right hand side of the estimate in Lemma \ref{dbar}, in view of Lemma \ref{benign} and the usual small constant - large constant estimate.

We first treat the case where $g\equiv 1$, in order to bring out more clearly the standard nature of the arguments involved. Let $u=\sum^{\prime}_{|J|=q}u_{J}d\overline{z_{J}}$ and
$v=\sum^{\prime}_{|K|=q+1}v_{K}d\overline{z_{K}}$.  We may assume that $X_{\epsilon} = e^{-h_{\epsilon}}\sum_{j=1}^{n}(\partial \rho/\partial
\overline{z_{j}}) (\partial/\partial z_{j})$ throughout $\Omega$. Indeed, the error is compactly supported and so is easily seen to be acceptable for the right hand side of \eqref{dbar0}, by interior elliptic regularity of $\overline{\partial} \oplus \overline{\partial}^{*}$. 
Then
\[ 
\left ( \left [\ \overline{\partial}, X_{\epsilon} \right ]u, v \right
)_{(0,q+1)} = \left ( \sum^{\prime}_{j,J}\left
(\sum_{k}\frac{\partial}{\partial \overline{z_{j}}}\left
(e^{-h_{\epsilon}}\frac{\partial \rho}{\partial \overline{z_{k}}}
\right )\frac{\partial u_{J}}{\partial z_{k}} \right )
d\overline{z_{j}} \wedge d\overline{z_{J}} \ , \ v \right )_{(0,q+1)}
\ .
\]
For $(j,J)$ fixed, the term $d\overline{z_{j}} \wedge
d\overline{z_{J}}$ picks out the component $v_{j,J}$ of $v$. So what
needs to be estimated is 
\begin{equation}\label{dbar1}
\sum^{\prime}_{j,J}\int_{\Omega}\left
(\sum_{k}\frac{\partial}{\partial
\overline{z_{j}}}\left (e^{-h_{\epsilon}}\frac{\partial
\rho}{\partial \overline{z_{k}}} \right )\frac{\partial
u_{J}}{\partial z_{k}} \right )\overline{v_{j,J}} \ .
\end{equation}
Note that $\partial u_{J}/\partial z_{k} = (\partial\rho)/(\partial
z_{k})\sum_{j=1}^{n}(\partial\rho/\partial 
\overline{z_{j}})(\partial u_{J}/ \partial z_{j}) + Y_{k}u_{J} = (\partial\rho/\partial z_{k})e^{h_{\epsilon}}X_{\epsilon}u_{J} +
Y_{k}u_{J}$ for a field $Y_{k}$ of type $(1,0)$ which is complex
tangential at the boundary and which does not depend on $\epsilon$.
The contribution coming from $Y_{k}u_{J}$ can be estimated using
Lemma \ref{benign}:
\begin{multline}\label{dbar1a}
\sum^{\prime}_{j,J}\int_{\Omega}\left
(\sum_{k}\frac{\partial}{\partial
\overline{z_{j}}}\left (e^{-h_{\epsilon}}\frac{\partial
\rho}{\partial \overline{z_{k}}} \right ) Y_{k}
u_{J}\right )\overline{v_{j,J}} \leq  C_{\epsilon}\sum_{k=1}^{n}
\|Y_{k}u\|\|v\| \\ 
 \leq C_{\epsilon}\left(
\|\overline{\partial}u\|^{2} + \|\overline{\partial}^{*}u\|^{2} +
\|u\|\|u\|_{1} \right )^{1/2}\|v\| \ .
\end{multline}
Using twice the inequality $|ab| \leq
(\delta/2)a^{2}+(1/2\delta)b^{2}$ together with $\|u\|^{2} \leq
C(\|\overline{\partial}u\|^{2}+\|\overline{\partial}^{*}u\|^{2})$ for
$u \in dom(\overline{\partial}^{*})$, the last expression is easily
seen to be bounded by the right hand side of the inequality in Lemma
\ref{dbar}. It remains to estimate the contribution in \eqref{dbar1}
that comes from the normal derivative of $u_{J}$. It equals
\begin{equation}\label{dbar2}
\sum^{\prime}_{j,J}\int_{\Omega}\left
(\sum_{k}\frac{\partial}{\partial
\overline{z_{j}}}\left (e^{-h_{\epsilon}}\frac{\partial
\rho}{\partial \overline{z_{k}}} \right ) \frac{\partial
\rho}{\partial z_{k}} e^{h_{\epsilon}}X_{\epsilon}
u_{J}\right )\overline{v_{j,J}} \ .
\end{equation}
Note that
$\sum_{k=1}^{n}(e^{-h_{\epsilon}}(\partial\rho/\partial\overline{z_{k}}))
((\partial\rho/\partial z_{k})e^{h_{\epsilon}}) \equiv 1/4$ near
$b\Omega$ (since $|\nabla\rho| \equiv 1$ near $b\Omega$). Therefore,
moving the derivative $\partial /\overline{\partial z_{j}}$ over from
one factor to the other in \eqref{dbar2} gives
\begin{equation}\label{dbar3}
-\sum_{J}^{\prime}\int_{\Omega}\left (
\sum_{j,k}\frac{\partial}{\partial\overline{z_{j}}}\left
(e^{h_{\epsilon}}\frac{\partial\rho}{\partial z_{k}}\right
)\frac{\partial\rho}{\partial\overline{z_{k}}}\;\overline{v_{jJ}}
\right )e^{-h_{\epsilon}}X_{\epsilon}u_{J} + O\left
(C_{\epsilon}\|v\|_{\Omega_{0}}\|X_{\epsilon}u\|_{\Omega_{0}} \right ) 
\end{equation}
for a suitable relatively compact subdomain $\Omega_{0}$. The compactly supported term can be estimated as at the beginning of the proof. Observe
that 
\[
\frac{\partial}{\partial\overline{z_{j}}}(e^{h_{\epsilon}}
\frac{\partial\rho}{\partial z_{k}}) = \frac{\partial^{2}}{\partial
z_{k}\partial\overline{z_{j}}}(e^{h_{\epsilon}}\rho) - \left
(\frac{\partial^{2}}{\partial
z_{k}\partial\overline{z_{j}}}e^{h_{\epsilon}}\right )\rho -
\frac{\partial}{\partial
z_{k}}(e^{h_{\epsilon}})\frac{\partial\rho}{\partial\overline{z_{j}}}
\ . 
\]
Inserting this into the first term in \eqref{dbar3} gives that this
term equals
\begin{multline}\label{dbar4}
-\sum_{J}^{\prime}\int_{\Omega}\left (
\sum_{j,k}\frac{\partial^{2}(e^{h_{\epsilon}}\rho)}
{\partial\overline{z_{j}}\partial z_{k}}
\;\frac{\partial\rho}{\partial\overline{z_{k}}}\;
\overline{v_{jJ}}\right )e^{-h_{\epsilon}}X_{\epsilon}u_{J} \\
+ \sum_{J}^{\prime}\int_{\Omega}\left (
\sum_{j,k}\frac{\partial^{2}e^{h_{\epsilon}}}
{\partial\overline{z_{j}}\partial z_{k}}
\;\frac{\partial\rho}{\partial\overline{z_{k}}}\;
\overline{v_{jJ}}\right )\rho e^{-h_{\epsilon}}X_{\epsilon}u_{J} \\
+ \sum_{J}^{\prime}\int_{\Omega}\left ( \overline{\left
(\sum_{j}\frac{\partial\rho}{\partial z_{j}}v_{j,J}\right )}
\left (\sum_{k}\frac{\partial\rho}{\partial
z_{k}}\frac{\partial}{\partial z_{k}}(e^{h_{\epsilon}})\right )
\right )e^{-h_{\epsilon}}X_{\epsilon}u_{J} \ .
\end{multline}
For the second term in \eqref{dbar4}, we have the upper bound
\begin{equation}\label{dbar5}
C_{\epsilon}\|v\|\|\rho X_{\epsilon}u\|  \leq
\epsilon\|v\|^{2} + C_{\epsilon}\|\rho X_{\epsilon}u\|^{2} \ .
\end{equation}
The first term on the right hand side of \eqref{dbar5} is dominated by
$\epsilon (\|\overline{\partial}v\|^{2} +
\|\overline{\partial}^{*}v\|^{2})$ (with a constant independent of
$\epsilon$). By an argument analogous to the
one used in estimating the second term in \eqref{hessian}, the second
term in \eqref{dbar5} is dominated by $\epsilon\|u\|_{1}^{2} +
C_{\epsilon} (\|\overline{\partial}u\|^{2} +
\|\overline{\partial}^{*}u\|^{2})$. Therefore, the right hand side of
\eqref{dbar5} is acceptable for the estimate in Lemma \ref{dbar}.

To estimate the third term in \eqref{dbar4}, first note that since
$v\in dom(\overline{\partial}^{*})$,
$\sum_{j}\frac{\partial\rho}{\partial z_{j}}v_{j,J} = 0$ on
$b\Omega$. Therefore,
\begin{multline}\label{dbar6}
\left\|\sum_{j}\frac{\partial\rho}{\partial z_{j}}v_{j,J}\right\|_{1}
\leq C 
\left\| \Delta
(\sum_{j}\frac{\partial\rho}{\partial z_{j}}v_{j,J})\right\|_{-1} \\
\leq 
C (\|\overline{\partial}v\| + \|\overline{\partial}^{*}v\| + \|v\|
+ \|v\|_{-1}) 
\leq C(\|\overline{\partial}v\| + \|\overline{\partial}^{*}v\|) \ .
\end{multline}
We have used in \eqref{dbar6} that
$\overline{\partial}\overline{\partial}^{*} +
\overline{\partial}^{*}\overline{\partial}$ acts componentwise as the
Laplacian (up to a constant), so that $\|\Delta v_{j,J}\|_{-1} \leq
\|\Delta v\|_{-1} \approx
\|(\overline{\partial}\overline{\partial}^{*} +
\overline{\partial}^{*}\overline{\partial})v\|_{-1} \leq
C(\|\overline{\partial}v\| + \|\overline{\partial}^{*}v\|)$. Using
the Cauchy-Schwarz inequality and interpolation of Sobolev norms, the
third term in \eqref{dbar4} can now be estimated by
\begin{multline}\label{dbar7}
C_{\epsilon}\sum_{J}^{\prime}\left\|\sum_{j}\frac{\partial\rho}{\partial
z_{j}}v_{j,J}\right\|\|u_{J}\|_{1} \leq \epsilon\|u\|_{1}^{2} +
C_{\epsilon}\sum_{J}^{\prime}\left\|\sum_{j}\frac{\partial\rho}{\partial
z_{j}}v_{j,J}\right\|^{2} \\
\leq \epsilon\|u\|_{1}^{2} +
\epsilon\sum_{J}^{\prime}\left\|\sum_{j}\frac{\partial\rho}{\partial
z_{j}}v_{j,J}\right\|_{1}^{2} +
C_{\epsilon}\sum_{J}^{\prime}\left\|\sum_{j}\frac{\partial\rho}{\partial
z_{j}}v_{j,J}\right\|_{-1}^{2} \\
\leq \epsilon\|u\|_{1}^{2} + \epsilon C(\|\overline{\partial}v\|^{2} +
\|\overline{\partial}^{*}v\|^{2}) + C_{\epsilon}\|v\|_{-1}^{2} \ .
\end{multline}
In the last step, we have used \eqref{dbar6}. Again, the right hand
side of \eqref{dbar7} is acceptable for the estimate in Lemma
\ref{dbar}.

It remains to consider the first term in \eqref{dbar4}. It is
estimated by 
\begin{multline}\label{dbar8}
\left\|\sum_{J}^{\prime}\left
(\sum_{j,k}\frac{\partial^{2}(e^{h_{\epsilon}}\rho)}
{\partial\overline{z_{j}}\partial z_{k}}
\;\frac{\partial\rho}{\partial\overline{z_{k}}}\;\overline{v_{jJ}}\right
) d\overline{z_{J}}\right\|
\left\|e^{-h_{\epsilon}}X_{\epsilon}u\right\| \\
\leq \sqrt{\epsilon}\|u\|_{1}^{2} +
(C/\sqrt{\epsilon})\left\|\sum_{J}^{\prime}\left
(\sum_{j,k}\frac{\partial^{2}(e^{h_{\epsilon}}\rho)}
{\partial\overline{z_{j}}\partial z_{k}}
\;\frac{\partial\rho}{\partial\overline{z_{k}}}\;\overline{v_{jJ}}\right
) d\overline{z_{J}}\right\|^{2} \ . 
\end{multline}
Note that by Lemma \ref{percol}, \eqref{ii} also holds for
$(q+1)$-forms. In particular, \eqref{ii} applies to $v$, that is, to
the last expression in \eqref{dbar8}, and this shows that the right
hand side of \eqref{dbar8} is bounded by the right hand side of the
estimate in Lemma \ref{dbar}. This concludes the proof of Lemma \ref{dbar} when $g \equiv 1$. For general $g$, keeping track of how $g$ enters the estimates, combined with standard arguments, gives the proof.  $\Box$
\smallskip

For commutators with $\overline{\partial}^{*}$, we let $X_{\epsilon}-\overline{X_{\epsilon}}$ act in special boundary charts, so that the domain of $\overline{\partial}^{*}$ is preserved. Denote by $\overline{\nabla}u$ the vector of all bar derivatives of all coefficients (say in Euclidean coordinates, although this is immaterial) of a form $u$.

\begin{lemma}\label{dbar*}
Let $1 \leq q \leq n$, and assume $\{\rho_{\epsilon}\}$ is a family of
defining functions as in Theorem \ref{main}, and let
$\{X_{\epsilon}\}$ be a family of smooth vector fields of type
$(1,0)$ so that near $b\Omega$ (possibly depending on $\epsilon$),
$X_{\epsilon}$ agrees with
$e^{-h_{\epsilon}}\sum_{j=1}^{n}(\partial \rho/\partial
\overline{z_{j}}) (\partial/\partial z_{j})$. Let $X_{\epsilon}-\overline{X_{\epsilon}}$ act componentwise in special boundary charts. There is a constant $A$, such that given a family of positive constants $B_{\epsilon}$, there are constants $C_{\epsilon, g, B_{\epsilon}}$, $0<\epsilon<1$, $g \in C^{\infty}(\overline{\Omega})$, such that for all $v \in
C_{(0,q)}^{\infty}(\overline{\Omega}) \cap
dom(\overline{\partial}^{*}), u \in
C_{(0,q-1)}^{\infty}(\overline{\Omega}) \cap
dom(\overline{\partial}^{*})$, we have the estimate 
\begin{multline}\label{dbar*1}
\left | \left (\left [\overline{\partial}^{*},X_{\epsilon}-\overline{X_{\epsilon}}
\right ](gv),u \right )_{(0,q-1)} \right | \\ 
\leq
A\sqrt{\epsilon}\|g\|_{\infty}^{2}\left
(\|\overline{\partial}v\|_{1}^{2}+ 
\|\overline{\partial}^{*}v\|_{1}^{2} +
\|v\|_{1}^{2} \right ) + 
A\sqrt{\epsilon} \left (\|u\|^{2}
+ \frac{1}{B_{\epsilon}}\|\overline{\nabla}u\|^{2} \right ) \\
+ C_{\epsilon, g, B_{\epsilon}} \left
(\|\overline{\partial}v\|^{2} + \|\overline{\partial}^{*}v\|^{2}\right) \ .
\end{multline}
\end{lemma}
\smallskip

\emph{Proof}: Let $v$ and $u$ as in the lemma. We will again use the standard estimate $\|v\|^{2} \leq C\left (\|\overline{\partial}v\|^{2} + \|\overline{\partial}^{*}v\|^{2}\right )$ for forms in $dom(\overline{\partial}) \cap dom(\overline{\partial}^{*})$ throughout the proof. Integration by parts gives
\begin{multline}\label{dbar*2}
\left | \left (\left [\overline{\partial}^{*},X_{\epsilon}-\overline{X_{\epsilon}}
\right ](gv),u \right )_{(0,q-1)} \right | \leq
\left | \left (gv,\left [\overline{\partial}^{\ }, X_{\epsilon}
- \overline{X_{\epsilon}} \right ]u 
\right )_{(0,q)}\right | \\
+ C_{\epsilon}\left ( \|gv\|\|\overline{\partial}u\| +
\|\overline{\partial}^{*}(gv)\|\|u\| \right ) \, .
\end{multline}
The last term on the right hand side of \eqref{dbar*2} is easily seen to be dominated by the right hand side of \eqref{dbar*1}. In the inner product on the right hand side of \eqref{dbar*2}, the contribution coming from $\overline{X_{\epsilon}}$ only involves bar derivatives of $u$, so is of order $C_{\epsilon}\|\overline{\nabla}u\|\|gv\| \leq
(\sqrt{\epsilon}/B_{\epsilon})\|\overline{\nabla} u\|^{2} + (B_{\epsilon}C_{\epsilon}^{2}/4\sqrt{\epsilon})\|g\|_{\infty}^{2}\|v\|^{2}$. This is dominated by the right hand side of \eqref{dbar*1}. To estimate the contribution from the commutator with $X_{\epsilon}$, we essentially repeat the proof of Lemma \ref{dbar}, but with the small constants - large constants arguments so that the norms involving $u$ appear only with small constants. Also, in several places derivatives will have to be integrated by parts to the other side of an inner product. Note that we may change the commutator to that with $X_{\epsilon}$ acting in Euclidean coordinates (which is the situation in Lemma \ref{dbar}): the error this makes is of order $C_{\epsilon}\|gv\|\left (\|\overline{\partial}u\|+\|u\| \right )$, which is acceptable for the right hand side of \eqref{dbar*1}. The details are as follows.  The tangential derivative $Y_{k}$ in \eqref{dbar1a} can be integrated by parts. The result is that this term is dominated by $C_{\epsilon}\|u\|\left (\|gv\| + \|Y_{k}(gv)\|\right )$. By Lemma \ref{benign}, this can be estimated by $\sqrt{\epsilon}\|u\|^{2} + \sqrt{\epsilon}\|g\|_{\infty}^{2}\|v\|_{1}^{2}+ C_{\epsilon, g}\left (\|\overline{\partial}v\|^{2} + \|\overline{\partial}^{*}v\|^{2}\right )$. We now proceed to \eqref{dbar2}; in turn, this leads to \eqref{dbar3}. After integrating $X_{\epsilon}$ by parts, we may write the estimate for the compactly supported term in \eqref{dbar3} as $C_{\epsilon}\left (\|\overline{X_{\epsilon}}(gv)\|_{\Omega_{0}} + \|gv\|_{\Omega_{0}}\right )\|u\|_{\Omega_{0}}$, which is acceptable for \eqref{dbar*1} (again by interior elliptic regularity). Proceeding to \eqref{dbar4}, we first consider the second term. Replacing $X_{\epsilon}$ by $X_{\epsilon}-\overline{X_{\epsilon}}$ makes an acceptable error of order $\|gv\| \| \overline{\nabla}u\|$. Integrating
$X_{\epsilon}-\overline{X_{\epsilon}}$ by parts shows that this term is bounded by 
\begin{equation}\label{dbar*3}
C_{\epsilon}\left (\|u\|\,\|gv\| + \|u\|\|\rho \left (X_{\epsilon}-\overline{X_{\epsilon}}\right )(gv)\|\right ) \, .
\end{equation}
The first term in \eqref{dbar*3} is acceptable. The second term is dominated by
\begin{multline}\label{dbar*4}
\epsilon \|u\|^{2} + C_{\epsilon}\|\rho \left (X_{\epsilon}-\overline{X_{\epsilon}}\right )(gv)\|^{2} \\
\leq
\epsilon \|u\|^{2} +
\epsilon\|g\|_{\infty}^{2}\|v\|_{1}^{2} +
C_{\epsilon, g}\left (\|\overline{\partial}v\|^{2} +
\|\overline{\partial}^{*}v\|^{2}\right ) \, ,
\end{multline}
and so is also acceptable. In the third term in \eqref{dbar4}, we again replace $X_{\epsilon}$ by $X_{\epsilon}-\overline{X_{\epsilon}}$, making an acceptable error. Integrating $ X_{\epsilon}-\overline{X_{\epsilon}} $ by parts gives a bound
\begin{equation}\label{dbar*5}
C_{\epsilon} \left (\|u\|\,\|gv\| + 
\|u\|\sum_{J}^{\prime}
\left \|\left (X_{\epsilon}-\overline{X_{\epsilon}}
\right ) \left (
g\sum_{j}\frac{\partial \rho}{\partial z_{j}}v_{j,J} \right )
\right \| \; \right ) \,.
\end{equation}
Applying \eqref{dbar6} to the second term on the right hand side of \eqref{dbar*5} shows that this right hand side can be bounded as required in \eqref{dbar*1}. To estimate the first term in \eqref{dbar4}, we once more replace $X_{\epsilon}$ by $X_{\epsilon}-\overline{X_{\epsilon}}$ (making an acceptable error) and integrate $X_{\epsilon}-\overline{X_{\epsilon}}$ by parts. The main term to be estimated is
\begin{multline}\label{dbar*6}
\left |\sum_{J}^{\prime}\int_{\Omega}
\left (\sum_{j,k}\frac{\partial^{2}\rho_{\epsilon}}{\partial 
\overline{z_{j}}\partial z_{k}}\frac{\partial \rho}{\partial \overline{z_{k}}}\overline{(X_{\epsilon}-\overline{X_{\epsilon}})(gv)_{j,J}} \right )u_{J} \right | \\
\leq
\left\|\sum_{J}^{\prime}
\left (\sum_{j,k}\frac{\partial^{2}\rho_{\epsilon}}{\partial 
\overline{z_{j}}\partial z_{k}}\frac{\partial \rho}{\partial \overline{z_{k}}}\overline{(X_{\epsilon}-\overline{X_{\epsilon}})(gv)_{j,J}} \right )d\overline{z_{J}}\right\|\|u\| \; .
\end{multline}
We may replace $(X_{\epsilon}-\overline{X_{\epsilon}})(gv_{j,J})$ in \eqref{dbar*6} by $T(gv_{j,J})$ where $T = L_{n} - \overline{L_{n}}$, and $L_{n} = \sum_{j=1}^{n} (\partial\rho/\partial\overline{z_{j}})(\partial/\partial z_{j})$. In order to apply \eqref{ii} to $T(gv)$, we need to switch back to $T$ acting in special boundary charts ( so that $T(gv) \in dom(\overline{\partial}^{*})$). We have $T(gv_{j,J}) =
(T(gv))_{j,J}$, where on the right $T$ acts on forms in special boundary charts, plus terms of order zero, which therefore are acceptable for \eqref{dbar*1}. Combining these observations and using \eqref{ii} (this time for $q$-forms) shows that the right hand side of \eqref{dbar*6} can be estimated by
\begin{multline}\label{dbar*7}
C_{\epsilon,g}\|u\|\,\|v\| + \frac{\sqrt{\epsilon}}{2}\|u\|^{2} +
\frac{1}{4\sqrt{\epsilon}}\left \|\sum_{J}^{\prime}
\left (\sum_{j,k}\frac{\partial^{2}\rho_{\epsilon}}
{\partial z_{k}\partial\overline{z_{j}}}
\frac{\partial \rho}{\partial \overline{z_{k}}}
\overline{(T(gv))_{j,J}}\right )d\overline{z_{J}} 
\right \|^{2}  \\
\leq
\sqrt{\epsilon}\|u\|^{2} + \frac{\sqrt{\epsilon}}{4}\|g\|_{\infty}^{2}
\left (\|\overline{\partial}(Tv)\|^{2}+
\|\overline{\partial}^{*}(Tv)\|^{2} \right )  +
C_{\epsilon,g}\left (\|Tv\|_{-1}^{2} + \|v\|^{2} \right )\\
\leq
\sqrt{\epsilon}\|u\|^{2} +
\frac{\sqrt{\epsilon}}{4}\|g\|_{\infty}^{2}
\left ( \|\overline{\partial}v\|_{1}^{2} +
\|\overline{\partial}^{*}v\|_{1}^{2} + 
\|v\|_{1}^{2} \right )  + C_{\epsilon,g}\|v\|^{2} \; .
\end{multline}
We have used here that $(T(gv))_{j,J} = (Tg)v_{j,J} + g(Tv)_{j,J}$ and that
$\left [\overline{\partial},T \right ]$ and $\left [\overline{\partial}^{*},T \right ]$ are operators of order one. This completes the proof of Lemma \ref{dbar*}. $\Box$

\section{Proof of Theorem \ref{main}}

We now come to the proof of Theorem \ref{main}. It uses many ideas from \cite{BoasStraube91}. In particular, we use a downward induction on the degree $q$. Fix a degree $q_{0}$, and assume that \eqref{ii} holds for $q_{0}$-forms, hence for $m$-forms for $m = q_{0}, \cdots, n$, by Lemma \ref{percol}. In top degree, $N_{n}$ is regular in Sobolev norms: the $\overline{\partial}$-Neumann boundary conditions reduce to Dirichlet boundary conditions, and the problem becomes coercive, see e.g. \cite{FollandKohn72},p.63. Therefore, to prove the theorem for $q_{0}$-forms, it suffices to show the following: if $N_{m}$ satisfies the Sobolev estimates \eqref{Sobolev} for $(q+1) \leq m \leq n$, and if \eqref{ii} holds for $q$-forms, then the estimates \eqref{Sobolev} hold for $q$-forms. We will use that Sobolev estimates for $N_{m}$ imply Sobolev estimates for the Bergman projection $P_{m-1}$ on $(m-1)$-forms (see \cite{BoasStraube90}); in particular, $P_{q}$ satisfies Sobolev estimates as a result of the induction assumption.

The arguments will involve absorbing terms, and one has to know that the terms to be absorbed are finite. Therefore, we first prove estimates for the regularized $\overline{\partial}$-Neumann operator $N_{\delta,q}$ for $\delta > 0$, where $N_{\delta,q}$ is the operator obtained from the usual elliptic regularization procedure (\cite{FollandKohn72}, section 3, chapter 2, \cite{Taylor96}, section 5, chapter 12). We will get the desired estimates for $N_{q}$ by letting $\delta$ tend to zero. Section 6 below contains various facts about the regularized problem that we will use. For the moment, we note that $N_{\delta,q}$ is the inverse of the selfadjoint operator $\Box_{\delta}$ associated to the quadratic form
\begin{equation}\label{qform}
Q_{\delta}(u,\overline{u}) = \|\overline{\partial}u\|^{2} +
\|\overline{\partial}^{*}u\|^{2} + \delta\|\nabla u\|^{2} 
\end{equation}
with form domain $W^{1}_{(0,q)}(\Omega) \cap dom(\overline{\partial}^{*})$, where $\nabla$ is the vector of all (first) derivatives of all components of $u$. As such, $N_{\delta,q}$ maps $\mathcal{L}^{2}_{(0,q)}(\Omega)$ continuously into this domain (endowed with the norm induced by $Q_{\delta}$).

We prove estimates for $\overline{\partial}N_{\delta,q}$, $\overline{\partial}^{*}N_{\delta,q}$, and $\delta^{1/2}\nabla N_{\delta, q}u$ (uniform in $\delta$ for small $\delta$) in $W^{k}(\Omega)$, by induction on $k$. That is, we will estimate $\|\overline{\partial}N_{\delta, q}u\|_{k}^{2} + \|\overline{\partial}^{*}N_{\delta, q}u\|_{k}^{2} + \delta\|\nabla N_{\delta, q}u\|_{k}^{2}$. By Lemma \ref{Reg2} in section 6, this also gives estimates, uniform in $\delta$, for $N_{\delta,q}$ in $W^{k}_{(0,q)}(\Omega)$. The case $k=0$ is taken care of by the preceeding remark. Let $u \in C^{\infty}_{(0,q)}(\overline{\Omega})$. Then $N_{\delta, q}u \in C^{\infty}_{(0,q)}(\overline{\Omega})$ for $\delta > 0$. As in \cite{BoasStraube91}, we use Lemma \ref{benign} to essentially reduce the problem to having to consider only tangential derivatives, and only in a direction transverse to the complex tangent space. The normal derivative of a form $u$ can be expressed in terms of tangential derivatives, components of $\overline{\partial}u$, of $\vartheta u$, and of $u$, where $\vartheta$ denotes the formal adjoint of $\overline{\partial}$ (the boundary is noncharacteristic for $\overline{\partial} \oplus \vartheta$).
Let $\rho_{\epsilon}$ be the family of defining functions given by the assumption in Theorem \ref{main} and choose a defining function $\rho$ with $|\nabla\rho| \equiv 1$ near $b\Omega$. Choose functions $h_{\epsilon}$ in $C^{\infty}(\overline{\Omega})$, bounded on $\Omega$ independently of $\epsilon$, so that  $\rho_{\epsilon}=e^{h_{\epsilon}}\rho$ near (depending on $\epsilon$) $b\Omega$. This is possible in view of \eqref{i}. Set $X_{\epsilon} := e^{-h_{\epsilon}}\sum_{j=1}^{n}(\partial\rho/\partial\overline{z_{j}})\partial/\partial z_{j}$. We let vector fields (derivatives) act on forms coefficientwise in special boundary charts. Then, tangential derivatives will preserve the domain of $\overline{\partial}^{*}$. Compactly supported terms can be bounded by interior elliptic regularity of $\,\Box_{\delta, q}$ (uniformly in $\delta$). Combining the previous remark with Lemma \ref{benign}, and absorbing $(s.c.)\|\overline{\partial}^{*}N_{\delta, q}u\|_{k}^{2}$, gives (compare also \cite{BoasStraube91}, p.83)
\begin{multline}\label{*N}
\|\overline{\partial}^{*}N_{\delta,q}u\|_{k}^{2} \leq C\left\|\left (X_{\epsilon}-\overline{X_{\epsilon}}\right )^{k} \overline{\partial}^{*}N_{\delta,q}u\right\|^{2} \\
+
C_{\epsilon}\left (\|\overline{\partial}\overline{\partial}^{*}N_{\delta,q}u\|_{k-1}^{2} +
\|\overline{\partial}^{*}N_{\delta,q}u\|_{k-1}^{2} + \|u\|_{k}^{2} \right ) \; .
\end{multline}
where $C$ does not depend on $\epsilon$ (because the $h_{\epsilon}$ are bounded on $\Omega$ independently of $\epsilon$). 

For $\overline{\partial}N_{\delta,q}$, the argument is more involved, because $\overline{\partial}N_{\delta,q}$ is not, in general, in the domain of $\overline{\partial}^{*}$ (so that Lemma \ref{benign} does not apply directly). Computing the free boundary condition for $\Box_{\delta,q}$ (see section 6 below) shows that the modified form $\overline{\partial}N_{\delta,q} + \delta
(\partial/\partial \nu)N_{\delta,q}u \wedge \overline{\omega_{n}}$ belongs to the domain of $\overline{\partial}^{*}$, where $\partial/\partial \nu$ denotes the normal derivative acting coefficientwise in Euclidean coordinates and $\omega_{n}$ is the $(1,0)$-form dual to $L_{n}=\sum_{j=1}^{n}(\partial\rho/\partial\overline{z_{j}})\partial/\partial z_{j}$. Applying the above reasoning to this modified form results in the estimate 
\begin{multline}\label{barN}
\left\|\overline{\partial}N_{\delta,q}u\right\|_{k}^{2} \leq C\left\|\left (X_{\epsilon}-\overline{X_{\epsilon}}\right )^{k} \overline{\partial}N_{\delta,q}u\right\|^{2}  \\
+
C_{\epsilon}\left (\left\|\vartheta  \overline{\partial}N_{\delta,q}u
\right\|_{k-1}^{2} +
\left\|\overline{\partial}N_{\delta,q}u
\right\|_{k-1}^{2} 
+ \left \|\delta (\partial / \partial \nu)N_{\delta,q}u \wedge \overline{\omega_{n}} \right\|_{k}^{2} + \|u\|_{k}^{2}     \right ) \; .
\end{multline}

For $\delta^{1/2}\nabla N_{\delta, q}u$, we obtain (via $\|\nabla N_{\delta, q}u\|_{k}^{2} \approx \|N_{\delta, q}u\|_{k+1}^{2}$)
\begin{multline}\label{nablaN}
\delta\|\nabla N_{\delta, q}u\|_{k}^{2} \leq C\delta\|(X_{\epsilon} - \overline{X_{\epsilon}})^{k}\nabla N_{\delta, q}u\|^{2} \\
+ C_{\epsilon}\delta\left ( \|\overline{\partial}N_{\delta, q}u\|_{k}^{2} + \|\overline{\partial}^{*}N_{\delta, q}u\|_{k}^{2} + \|N_{\delta, q}u\|_{k}\|N_{\delta, q}u\|_{k+1} + \|u\|_{k}^{2} \right ) \; .
\end{multline}

The terms $\|\overline{\partial}N_{\delta, q}u\|_{k-1}^{2}$ and $\|\overline{\partial}^{*}N_{\delta, q}u\|_{k-1}^{2}$ in \eqref{*N} and \eqref{barN} are bounded by $\|u\|_{k-1}^{2}$, by induction assumption. The second to the last term in \eqref{barN} is of order $C_{\epsilon}\delta^{2}\|\nabla N_{\delta, q}u\|_{k}^{2}$. Upon adding \eqref{*N}, \eqref{barN}, and \eqref{nablaN}, it can be absorbed for $\delta < \delta(\epsilon)$ (there is an extra factor $\delta$). By Lemma \ref{Reg4}, we have that $\|\overline{\partial}\overline{\partial}^{*}N_{\delta, q}u\|_{k-1}^{2}  +   \|\vartheta\overline{\partial}N_{\delta, q}u\|_{k-1}^{2}$ is dominated by $\|\overline{\partial}N_{\delta, q}u\|_{k-1}^{2}   +  \delta\|\nabla N_{\delta, q}u\|_{k-1}^{2}  +  \|u\|_{k-1}^{2} + \|u\|_{1}^{2} + \delta^{2}\|u\|_{k}^{2}$, which in turn is dominated by $\|u\|_{k-1}^{2} + \|u\|_{1}^{2} + \delta^{2}\|u\|_{k}^{2} \leq 3\|u\|_{k}^{2}$, by induction assumption. The terms in \eqref{nablaN} involving $\|\overline{\partial}N_{\delta, q}u\|_{k}^{2}$ and $\|\overline{\partial}^{*}N_{\delta, q}u\|_{k}^{2}$ can be absorbed for $\delta < \delta(\epsilon)$. Finally, we note that $C_{\epsilon}\delta\|N_{\delta, q}u\|_{k}\|N_{\delta, q}u\|_{k+1} \leq (C_{\epsilon}^{2}/2)\delta\|N_{\delta, q}u\|_{k}^{2}  +  (1/2)\delta\|\nabla N_{\delta, q}u\|_{k}^{2}$. This, in view of Lemma \ref{Reg2}, can be absorbed into the sum of the left hand sides of \eqref{*N}, \eqref{barN}, and \eqref{nablaN}, again for $\delta < \delta(\epsilon)$. Thus, what remains to be estimated is $\| (X_{\epsilon}-\overline{X_{\epsilon}})^{k} \overline{\partial}N_{\delta,q}u \|^{2} +
\| (X_{\epsilon}-\overline{X_{\epsilon}} )^{k} \overline{\partial}^{*}N_{\delta,q}u \|^{2} +
\delta\|(X_{\epsilon} - \overline{X_{\epsilon}})^{k}\nabla N_{\delta, q}u\|^{2}$. 

We have
\begin{multline}\label{est1}
\left ( \left (X_{\epsilon}-\overline{X_{\epsilon}}\right )^{k}
\overline{\partial}^{*}N_{\delta,q}u, \left (X_{\epsilon}-\overline{X_{\epsilon}}\right )^{k}
\overline{\partial}^{*}N_{\delta,q}u \right ) \\
= \left (\overline{\partial}^{*}N_{\delta,q}u, 
\left (X_{\epsilon}-\overline{X_{\epsilon}}\right )^{2k}
\overline{\partial}^{*}N_{\delta,q}u \right ) +
O_{\epsilon}\left (\|\overline{\partial}^{*}N_{\delta,q}u\|_{k-1}
\|\overline{\partial}^{*}N_{\delta,q}u\|_{k} \right ) \; .
\end{multline}
The first term on the right hand side in \eqref{est1} equals
\begin{multline}\label{est2}
\left (\overline{\partial}^{*}N_{\delta,q}u, 
\overline{\partial}^{*}\left (X_{\epsilon}-\overline{X_{\epsilon}}\right )^{2k}
N_{\delta,q}u \right ) + 
\left (\overline{\partial}^{*}N_{\delta,q}u, 
\left [ \left (X_{\epsilon}-\overline{X_{\epsilon}}\right )^{2k},
\overline{\partial}^{*} \right ]N_{\delta,q}u \right ) \; .
\end{multline}
Expanding the commutator $\left [ \left (X_{\epsilon}-\overline{X_{\epsilon}}\right )^{2k},
\overline{\partial}^{*} \right ]$
in the usual way (see e.g. \cite{DerridjTartakoff76}, Lemma 2, p.418) gives
\begin{multline}\label{est3}
\left ( \left (X_{\epsilon}-\overline{X_{\epsilon}}\right )^{k}
\overline{\partial}^{*}N_{\delta,q}u, \left (X_{\epsilon}-\overline{X_{\epsilon}}\right )^{k}
\overline{\partial}^{*}N_{\delta,q}u \right ) \\
=\left (\overline{\partial}^{*}N_{\delta,q}u, 
\overline{\partial}^{*}\left(X_{\epsilon}-\overline{X_{\epsilon}}\right )^{2k}
N_{\delta,q}u \right ) \ \ \ \ \ \ \ \ \  \ \ \ \ \ \ \ \ \ \ \ \ \ \ \ \ \ \ \ \ \ \ \ \ 
\ \ \ \ \ \ \ \ \ \ \ \ \ \ \  \\
\ \ \ \  \  + \left ( \left (X_{\epsilon}-\overline{X_{\epsilon}}\right )^{k}
\overline{\partial}^{*}N_{\delta,q}u,\left [ \left (X_{\epsilon}-\overline{X_{\epsilon}}\right ),
\overline{\partial}^{*} \right ]\left (X_{\epsilon}-\overline{X_{\epsilon}}\right )^{k-1}N_{\delta,q}u
\right ) \\
+ O_{\epsilon}\left (\|\overline{\partial}^{*}N_{\delta,q}u\|_{k-1}
\|\overline{\partial}^{*}N_{\delta,q}u\|_{k} \right ) \; .
\end{multline}
Note that we can always integrate powers of $(X_{\epsilon}-\overline{X_{\epsilon}})$ by parts back to the left hand side before applying the Cauchy-Schwarz inequality.

A similar computation for $\left ( \left (X_{\epsilon}-\overline{X_{\epsilon}}\right )^{k}
\overline{\partial}N_{\delta,q}u, \left (X_{\epsilon}-\overline{X_{\epsilon}}\right )^{k}
\overline{\partial}N_{\delta,q}u \right )$ gives
\begin{multline}\label{est4}
\left ( \left (X_{\epsilon}-\overline{X_{\epsilon}}\right )^{k}
\overline{\partial}N_{\delta,q}u, \left (X_{\epsilon}-\overline{X_{\epsilon}}\right )^{k}
\overline{\partial}N_{\delta,q}u \right ) \\
= \left (\overline{\partial}N_{\delta,q}u, 
\overline{\partial}\left (X_{\epsilon}-\overline{X_{\epsilon}}\right )^{2k}
N_{\delta,q}u \right ) \ \ \ \ \ \ \ \ \ \ \ \ \ \ \ \ \ \ \ \ \ \ \ \ \ \ \ \ \ \ \ \ \ \ \ \ \ \ \ \ \ \ \ \ \ \ \ \ \ \\
\ \ +
\left ( \left (X_{\epsilon}-\overline{X_{\epsilon}}\right )^{k}
\overline{\partial}N_{\delta,q}u,\left [ \left (X_{\epsilon}-\overline{X_{\epsilon}}\right ),
\overline{\partial} \right ]\left (X_{\epsilon}-\overline{X_{\epsilon}}\right )^{k-1}N_{\delta,q}u
\right ) \\
+ O_{\epsilon}\left (\left\|\overline{\partial}N_{\delta,q}u\right\|_{k-1}
\left\|\overline{\partial}N_{\delta,q}u\right\|_{k} \right ) \; .
\end{multline}

Likewise,
\begin{multline}\label{est4a}
\delta\left ((X_{\epsilon} - \overline{X_{\epsilon}})^{k}\nabla N_{\delta, q}u, (X_{\epsilon} - \overline{X_{\epsilon}})^{k}\nabla N_{\delta, q}u \right ) \\
= \delta\left (\nabla N_{\delta, q}u, \nabla (X_{\epsilon} - \overline{X_{\epsilon}})^{2k}N_{\delta, q}u \right ) \ \ \ \ \ \ \ \ \ \ \ \ \ \ \ \ \ \ \ \ \ \ \ \ \ \ \ \ \ \ \ \ \ \ \ \ \ \ \ \ \ \ \ \\
+
\delta\left ((X_{\epsilon}-\overline{X_{\epsilon}})^{k}\nabla N_{\delta, q}u,
\left [X_{\epsilon}-\overline{X_{\epsilon}}, \nabla\right ](X_{\epsilon}-\overline{X_{\epsilon}})^{k-1}N_{\delta, q}u \right ) \\
+ \delta O_{\epsilon}\left ( \|\nabla N_{\delta, q}u\|_{k-1}\|\nabla N_{\delta, q}u\|_{k}\right ) \; .
\end{multline}

When we add \eqref{*N}, \eqref{barN}, and \eqref{nablaN} and use estimates \eqref{est3}, \eqref{est4}, and \eqref{est4a}, the first terms on the right hand sides of these estimates  add up to $Q_{\delta}(N_{\delta,q}u, (X_{\epsilon}-\overline{X_{\epsilon}})^{2k}N_{\delta, q}u) = (u, (X_{\epsilon}-\overline{X_{\epsilon}})^{2k}N_{\delta, q}u)$; and their sum is dominated by $C_{\epsilon}\|u\|_{k}\|N_{\delta, q}u\|_{k} \leq C_{\epsilon}\|u\|_{k}(\|\overline{\partial}N_{\delta, q}u\|_{k} + \|\overline{\partial}^{*}N_{\delta, q}u\|_{k})$, again by Lemma \ref{Reg2}. Thus this term can be absorbed.

It remains to estimate the inner products on the right hand sides of \eqref{est3}, \eqref{est4}, and \eqref{est4a} which involve commutators with $\overline{\partial}^{*}$, $\overline{\partial}$, and $\nabla$, respectively. We begin with \eqref{est3}. Observe that $\left (X_{\epsilon}-\overline{X_{\epsilon}}\right )^{k-1}N_{\delta, q}u = e^{-(k-1)h_{\epsilon}}T^{k-1}N_{\delta, q}u + D_{\epsilon}^{k-2}N_{\delta, q}u$, where $T=L_{n} - \overline{L_{n}}$, and $D_{\epsilon}^{k-2}$ denotes a differential operator of order $(k-2)$ with coefficients depending on $\epsilon$ (more precisely, on $h_{\epsilon}$ and its derivatives up to order $k-1$). The contribution coming from this term is $O_{\epsilon} \left ( \|\overline{\partial}^{*}N_{\delta, q}u  \|_{k}\,\|N_{\delta, q}u  \|_{k-1} \right ) \leq s.c.  \|\overline{\partial}^{*}N_{\delta, q}u  \|_{k}^{2} + (l.c.)_{\epsilon} \|N_{\delta, q}u \|_{k-1}^{2}$. The first term on the right can be absorbed, the second is dominated by $\|\overline{\partial}N_{\delta, q}u\|_{k-1}^{2} + \|\overline{\partial}^{*}N_{\delta, q}u\|_{k-1}^{2}$ (as above) and thus by $\|u\|_{k-1}^{2}$, by induction assumption.. To estimate the contribution from $e^{-(k-1)h_{\epsilon}}T^{k-1}N_{\delta, q}u$, we apply Lemma \ref{dbar*} with $g=e^{-(k-1)h_{\epsilon}}$ and with the family $B_{\epsilon}$ to be specified below. Then $\|g\|_{\infty} = \|e^{-(k-1)h_{\epsilon}}\|_{\infty} \leq C $ independently of $\epsilon$ ($k$ is fixed). Thus this contribution is dominated by
\begin{multline}\label{est7}
\sqrt{\epsilon} \left (  \|\overline{\partial}T^{k-1}N_{\delta, q}u \|_{1}^{2}
+ \|\overline{\partial}^{*}T^{k-1}N_{\delta, q}u \|_{1}^{2}
+ \|T^{k-1}N_{\delta, q}u \|_{1}^{2} \right.  \\
\left. + \| \left (X_{\epsilon}-\overline{X_{\epsilon}} \right )^{k} \overline{\partial}^{*}N_{\delta, q}u \|^{2}
+ (1 / B_{\epsilon})\|\overline{\nabla}\left (X_{\epsilon}-\overline{X_{\epsilon}} \right )^{k}\overline{\partial}^{*}N_{\delta, q}u\|^{2} \right ) \\
+ C_{\epsilon, B_{\epsilon}} \left (\|\overline{\partial}T^{k-1}N_{\delta, q}u \|^{2}
+ \|\overline{\partial}^{*}T^{k-1}N_{\delta, q}u\|^{2} \right ) \; .
\end{multline}
Commuting $T$ with $\overline{\partial}$ and with $\overline{\partial}^{*}$, respectively, shows that the terms in the first line of \eqref{est7} are of order $\sqrt{\epsilon}(\|\overline{\partial}N_{\delta, q}u\|_{k}^{2} +
\|\overline{\partial}^{*}N_{\delta, q}u\|_{k}^{2} + \|N_{\delta, q}u\|_{k}^{2} )$ and so can be absorbed (again upon adding \eqref{*N}, \eqref{barN}, and \eqref{nablaN}, and for $\epsilon$ small enough). The first term in the second line of \eqref{est7} can be absorbed into the left hand side of \eqref{est3}. The terms in the third line of \eqref{est7}, again upon commuting $\overline{\partial}$ and $\overline{\partial}^{*}$ with $T$, are of lower order, and are handled by the induction assumption as above. Finally, the second term on the second line in \eqref{est7} dictates the choice of $B_{\epsilon}$: we choose $B_{\epsilon}$ big enough so that it dominates the coefficients of both the commutator $\left [\overline{\nabla}, (X_{\epsilon} - \overline{X_{\epsilon}})^{k} \right ]$ and $(X_{\epsilon} - \overline{X_{\epsilon}})^{k}$. Then, by Lemma \ref{benign}, this term is dominated, independently of $\epsilon$, by $\sqrt{\epsilon}(\|\overline{\partial}^{*}N_{\delta, q}u\|_{k}^{2} +
\|\overline{\partial}\overline{\partial}^{*}N_{\delta, q}u\|_{k}^{2}) \leq 
C\sqrt{\epsilon}(\|\overline{\partial}^{*}N_{\delta, q}u\|_{k}^{2} +
\|\overline{\partial}N_{\delta, q}u\|_{k}^{2} + \delta\|\nabla N_{\delta, q}u\|_{k}^{2} + \|u\|_{k}^{2} + \delta^{2}\|u\|_{k+1}^{2} )$, in view of Lemma \ref{Reg4} (note that here $k \geq 1$). The first three terms can be absorbed when $\epsilon$ is chosen small enough.

We now come to the term in \eqref{est4} that contains the commutator $[(X_{\epsilon} - \overline{X_{\epsilon}}), \overline{\partial}]$. We replace $\overline{\partial}N_{\delta, q}u$ by $\overline{\partial}N_{\delta, q}u + \delta(\partial /\partial\nu)N_{\delta, q}u \wedge \overline{\omega_{n}}$ (so as to be in $dom(\overline{\partial}^{*})$), making an error that is of order $O_{\epsilon}(\|\delta(\partial /\partial\nu)N_{\delta, q}u\|_{k} \|N_{\delta, q}u\|_{k}) \leq C_{\epsilon}(\delta^{1/2}\delta
\|\nabla N_{\delta, q}u\|_{k}^{2} + \delta^{1/2}\|N_{\delta, q}u\|_{k}^{2})$. The first term can be absorbed if $\delta < \delta(\epsilon)$; so can the second in view of Lemma \ref{Reg2}. Now we apply Lemma \ref{dbar}, and by arguments similar to the ones just completed, we obtain that the term is dominated by
\begin{multline}\label{est10}
\sqrt{\epsilon}\left (\|\overline{\partial}N_{\delta, q}u\|_{k}^{2} + 
\|\vartheta\overline{\partial}N_{\delta, q}u\|_{k}^{2} +
\|N_{\delta, q}u\|_{k}^{2} + \delta^{2}\|(\partial/\partial\nu)N_{\delta, q}u\|_{k+1}^{2}
\right )  \\ 
+ C_{\epsilon}\left (\|\overline{\partial}N_{\delta, q}u\|_{k-1}^{2}
+ \|\overline{\partial}^{*}N_{\delta, q}u\|_{k-1}^{2}
+ \|N_{\delta, q}u\|_{k-1}^{2}
+ \|\delta\nabla N_{\delta, q}u\|_{k-1}^{2} \right ) \; .
\end{multline}
The terms that come with the factor $C_{\epsilon}$ can be handled by the induction hypothesis on Sobolev estimates in $W^{k-1}$ (again in conjunction with Lemma \ref{Reg2} for the term $\|N_{\delta, q}u\|_{k-1}^{2}$). The terms that come with a factor $\sqrt{\epsilon}$ are bounded by $\|\overline{\partial}N_{\delta, q}u\|_{k}^{2} +
\|\overline{\partial}^{*}N_{\delta, q}u\|_{k}^{2} +
\delta\|\nabla N_{\delta, q}u\|_{k}^{2} + \|u\|_{k}^{2} + \delta^{2}\|u\|_{k+1}^{2}$, after applying Lemmas \ref{Reg2}, \ref{Reg3}, and \ref{Reg4}. Again, the terms that matter can be absorbed when $\epsilon$ is small enough, while we keep the terms that involve norms of $u$.

The term in \eqref{est4a} that contains the commutator $[X_{\epsilon}-\overline{X_{\epsilon}}, \nabla ]$ is easily seen to be dominated by 
$C_{\epsilon}\delta\|\nabla N_{\delta, q}u\|_{k}\|N_{\delta, q}u\|_{k} \leq
C_{\epsilon}\delta((s.c.)\|\nabla N_{\delta, q}u\|_{k}^{2} +
(l.c.)\|N_{\delta, q}u\|_{k}^{2})$. By choosing $(s.c.)$ small enough (depending on $\epsilon$), we can absorb the first term into the sum of the left hand sides of \eqref{*N}, \eqref{barN}, and \eqref{nablaN}. The second term, upon applying Lemma \ref{Reg2}, can be absorbed for $\delta < \delta(\epsilon)$.

Adding \eqref{*N}, \eqref{barN}, and \eqref{nablaN}, using the induction assumption on estimates in $W^{k-1}$, choosing $\epsilon > 0$ small enough and absorbing terms, we find
\begin{equation}\label{est11}
\|\overline{\partial}N_{\delta, q}u\|_{k}^{2}
+ \|\overline{\partial}^{*}N_{\delta, q}u\|_{k}^{2} 
+ \delta\|\nabla N_{\delta, q}u\|_{k}^{2} \leq C(\|u\|_{k}^{2} + \delta^{2}\|u\|_{k+1}^{2}) \,,\; \delta \leq \delta_{0} \; ,
\end{equation}
where $C$ is independent of $\delta$, and $\delta_{0} = \delta(\epsilon)$ is determined now that $\epsilon$ has been chosen. By Lemma \ref{Reg2}, \eqref{est11} implies the same estimate for $N_{\delta, q}u$, namely $\|N_{\delta, q}u\|_{k}^{2} \leq C(\|u\|_{k}^{2} + \delta^{2}\|u\|_{k+1}^{2}) \,, \; \delta \leq \delta_{0} \;$ , again with a constant $C$ that is independent of $\delta$. Letting $\delta \rightarrow 0^{+}$ gives the estimate $\|N_{q}u\|_{k} \leq C\|u\|_{k}$. Indeed, a subsequence of $N_{\delta, q}u$ converges weakly in $W_{(0,q)}^{k}(\Omega)$. This weak limit equals $N_{q}u$. This follows from the identity $(u,v) = Q_{\delta}(N_{\delta, q}u,v) = 
Q(N_{q}u,v)$ for $u,v \in C_{(0,q)}^{\infty}(\overline{\Omega}) \cap dom(\overline{\partial}^{*})$; compare \cite{ChenShaw01}, p.103. Therefore, $\|N_{q}u\|_{k}^{2} \leq \limsup_{\delta \rightarrow 0^{+}}C(\|u\|_{k}^{2} + \delta^{2}\|u\|_{k+1}^{2}) = C\|u\|_{k}^{2}$. The Sobolev estimate we have shown is for $u \in C^{\infty}_{(0,q)}(\overline{\Omega})$, but this space is dense in $W^{k}_{(0,q)}(\Omega)$, and $N_{q}$ is continuous in $\mathcal{L}^{2}_{(0,q)}(\Omega)$, so that the estimate carries over to $u \in W^{k}_{(0,q)}(\Omega)$. This completes the downward induction step from $(q+1)$ to $q$, and thus the proof of Theorem \ref{main}. $\Box$
\smallskip

The proof of Theorem \ref{main} is more closely inspired by \cite{BoasStraube91} than meets the eye. In fact, it is fairly easy to combine the arguments in \cite{BoasStraube91} with those in the proof of Lemma \ref{dbar} to obtain a much shorter proof, but only of a priori estimates. It is in turning these a priori estimates into genuine estimates that difficulties arise, quite in contrast to \cite{BoasStraube91}. The simple method of passing to interior approximating strictly pseudoconvex domains employed in \cite{BoasStraube91} does not seem to be applicable in our situation, as it is not clear whether our weaker assumptions are inherited by such subdomains. To accommodate elliptic regularization, it seemed advantageous to somewhat rearrange the arguments; the result is the above proof.

\section{Vector fields that commute approximately with $\overline{\partial}$}

In \cite{StraubeSucheston02}, the authors showed that several conditions, known to be sufficient for global regularity of the $\overline{\partial}$-Neumann problem, can be modified in a natural way so as to become equivalent (and still imply global regularity). The purpose of this section is to show that the (equivalent) modified conditions imply the condition in Theorem \ref{main}; that is, Theorem \ref{main} covers this approach to global regularity as well.

We recall two of the definitions from \cite{StraubeSucheston02}. Denote by $K$ the set of boundary points of infinite type (in the sense of D'Angelo, \cite{D'Angelo93}) of a smooth bounded pseudoconvex domain $\Omega$. Then $K$ is compact.

We say that $\Omega$ \emph{admits a family of vector fields transverse to $b\Omega$ that commutes approximately with $\overline{\partial}$ at points of $K$} if the following holds. There exists a constant $C > 0$ such that for every $\epsilon > 0$, there exists a vector field $X_{\epsilon}$ of type $(1,0)$ whose coefficients are smooth in a neighborhood (in $\mathbb{C}^{n}$) $U_{\epsilon}$ of $K$ and such that 
\begin{equation}\label{field1}
C^{-1} < |(X_{\epsilon}\rho)(z)| < C, \;\; |arg((X_{\epsilon}\rho)(z))| < \epsilon \;\;, z \in K \; ,
\end{equation}
and
\begin{equation}\label{field2}
\left |\partial\rho\left (\left [X_{\epsilon}, \partial /\partial \overline{z_{j}}\right ]
\right )(z)\right | < \epsilon\;, z \in K\;, 1 \leq j \leq n\;.
\end{equation}
Actually this is precisely the condition in the `vector field method' in \cite{BoasStraube93}. In \cite{StraubeSucheston02}, the sightly more restrictive definition was used that $X_{\epsilon}\rho(z)$ should be real when $z \in K$. We will indicate below that this is irrelevant.

The existence of the vector fields need not imply a defining function whose Hessian is positive semi-definite at boundary points, see remark 3, p.234, in \cite{BoasStraube93}. However, it does imply a defining function with a weaker property, which is still sufficient for global regularity. We have the following definition from \cite{StraubeSucheston02}. $\Omega$ \emph{admits a family of essentially pluriharmonic defining functions} if there exists $C>0$ such that for all $\epsilon > 0$ there is a defining function $\rho_{\epsilon}$ for $\Omega$ satisfying
\begin{equation}\label{def1}
C^{-1} \leq |\nabla\rho_{\epsilon}(z)| \leq C\;, \;z \in b\Omega \;,
\end{equation}
and
\begin{equation}\label{def2}
\left |\sum_{j,k}\frac{\partial^{2}\rho_{\epsilon}(P)}{\partial z_{j} \partial \overline{z_{k}}}w_{j}\overline{w_{k}} \right | \leq O(\epsilon)|w|^{2}\;, 
\forall w \in span_{\mathbb{C}}\{N(P), L_{n}(P)\}
\end{equation}
for all boundary points $P$ in $K$. $Span_{\mathbb{C}}$ denotes the linear span over $\mathbb{C}$, $N(P)$ is the nullspace of the Levi form at $P \in b\Omega$, 
and $L_{n} = \sum_{j=1}^{n}(\partial\rho /\partial\overline{z_{j}})(\partial /\partial z_{j})$ (for a fixed defining function $\rho$). We emphasize that this condition is indeed a generalization of the notion of a defining function that is pluri\emph{sub}harmonic at boundary points. That is, if a domain admits a defining function whose complex Hessian is positive semi-definite in all directions at points of the boundary, then it admits a family of essentially pluriharmonic defining functions; this is explained in detail in \cite{StraubeSucheston02}, p. 251-252, to where we refer the reader. 

The main result in \cite{StraubeSucheston02} says that $\Omega$ admits a family of vector fields transverse to $b\Omega$ that commutes approximately with $\overline{\partial}$ at points of $K$ if and only if $\Omega$ admits a family of essentially pluriharmonic defining functions (with the slightly more stringent definition pointed out above; we will take care of this point below). More is done in \cite{StraubeSucheston02}: these two properties are also equivalent to a suitably formulated approximate exactness property of the winding form (the form $\alpha$ in \cite{BoasStraube93}) in weakly pseudoconvex directions, as well as to the existence of normals which are approximately conjugate holomorphic in weakly pseudoconvex directions. We do not discuss this here and refer the reader to \cite{StraubeSucheston02}. 

We can now formulate the main result of this section.
\begin{proposition}\label{vectorfields}
Let $\Omega$ be a smooth bounded pseudoconvex domain in $\mathbb{C}^{n}$, denote by $K$ the set of boundary points of $\Omega$ of infinite type. Assume that $b\Omega$ admits a family of vector fields transverse to $b\Omega$ that commutes approximately with $\overline{\partial}$ at points of $K$. Then the conditions in Theorem \ref{main} are satisfied for $q=1, 2, \cdots , n$.
\end{proposition}

\emph{Proof:} By Lemma \ref{percol}, we only have to consider the case $q=1$. We use the result from \cite{StraubeSucheston02} that under the assumption in Proposition \ref{vectorfields}, $\Omega$ admits a family of essentially pluriharmonic defining functions, say $\{\rho_{\epsilon}\}_{\epsilon > 0}$. Polarization gives that
\begin{equation}\label{polarization}
\left |\sum_{j,k = 1}^{n} \frac{\partial^{2}\rho_{\epsilon}}{\partial z_{j} \partial\overline{z_{k}}}(P)\frac{\partial\rho}{\partial\overline{z_{j}}}(P)
\overline{w_{k}}
\right |^{2} \leq C\epsilon|w|^{2}\;, \forall w \in span_{\mathbb{C}}\{N(P), L_{n}(P)\}\;,
P \in K\; .
\end{equation}
\eqref{polarization} holds in particular for $(P,w) \in \widetilde{K}$, where $\widetilde{K}$ denotes the compact set $\{(P,w) /\; P \in K, |w| = 1, w \in N(P) \}$, viewed as a subset of the unit sphere bundle in the complex tangent space bundle to $b\Omega$. Choose an open neighborhood $V_{\epsilon}$ of $\widetilde{K}$ in this bundle such that \eqref{polarization} still holds ( with a bigger constant) when $(P,w) \in V_{\epsilon}$. There is a constant $C_{\epsilon}$ such that $|w|^{2} \leq C_{\epsilon}\sum_{j,k =1}^{n}(\partial^{2}\rho_{\epsilon} / (\partial z_{j}\partial\overline{z_{k}})(P)w_{j}\overline{w_{k}}$ when $P \in K, w \in T_{P}^{\mathbb{C}}(b\Omega), (P,w/|w|) \notin V_{\epsilon}$. Consequently, we have the following estimate when $P \in K, w \in T_{P}^{\mathbb{C}}(b\Omega)$:
\begin{equation}\label{mixed}
\left |\sum_{j.k =1}^{n}\frac{\partial^{2}\rho_{\epsilon}}{\partial z_{j}\partial\overline{z_{k}}}(P)\frac{\partial\rho}{\partial\overline{z_{j}}}(P)\overline{w_{k}} \right |^{2} \leq C\epsilon |w|^{2} + \widetilde{C_{\epsilon}}\sum_{j.k =1}^{n}\frac{\partial^{2}\rho_{\epsilon}}{\partial z_{j}\partial\overline{z_{k}}}(P)w_{j}\overline{w_{k}} \; .
\end{equation}
(Both terms on the right hand side of \eqref{mixed} are nonnegative; the first term dominates the left hand side when $(P,w/|w|) \in V_{\epsilon}$, the second term dominates $|w|^{2}$, hence the left hand side of \eqref{mixed}, when $(P,w/|w|) \notin V_{\epsilon}$; when $w=0$, there is nothing to prove.) By continuity and homogeneity, \eqref{mixed} holds (up to increasing $\epsilon$ to $2\epsilon$, etc.) for $z$ in a neighborhood (in $\mathbb{C}^{n}$) $W_{\epsilon}$ of $K$. Now let $u \in C^{\infty}_{(0,1)}(\overline{\Omega})$. Choose a smooth cutoff function $\varphi_{\epsilon}$ that is identically $1$ on $K$ and is supported in $W_{\epsilon}$. Then the contribution to the left hand side of \eqref{ii} coming from $(1-\varphi_{\epsilon})u$ can be dominated as required by the right hand side of \eqref{ii} by using subelliptic estimates on the support of $(1-\varphi_{\epsilon})$ and interpolation of Sobolev norms (compare the argument in \eqref{dbar7} above). The contribution from $\varphi_{\epsilon}u$ can be dominated by (note that $|\varphi_{\epsilon}| \leq 1$)
\begin{equation}\label{fields3}
C\epsilon\|u\|^{2} + \widetilde{C_{\epsilon}}\int_{\Omega}
\sum_{j,k =1}^{n}\frac{\partial^{2}\rho_{\epsilon}}{\partial z_{j}\partial\overline{z_{k}}}(z)u_{j}(z)\overline{u_{k}(z)} \; .
\end{equation}
By the discussion in section 2 (see in particular \eqref{hessian}), \eqref{fields3} can be bounded in the way required in Theorem \ref{main}. This completes the proof of Proposition \ref{vectorfields}. $\Box$

We now take the opportunity to clarify a point left open in \cite{StraubeSucheston03}: if $\Omega$ admits a family of vector fields transverse to the boundary that commutes approximately with $\overline{\partial}$ at points of $K$, then $\Omega$ admits a family of essentially pluriharmonic defining functions. This was shown in \cite{StraubeSucheston03} only under the slightly stronger assumption that $X_{\epsilon}\rho$ is real on $K$ (rather than only approximately real). Actually everything that is needed is in place in \cite{StraubeSucheston03}. Namely, the proof of the implication $(iii) \rightarrow (i)$ in the Theorem in \cite{StraubeSucheston02}, when followed verbatim, gives a certain `defining function' $\widehat{\rho_{\epsilon}} = e^{h_{\epsilon}}\rho$. Here, $h_{\epsilon}$ is defined on $b\Omega$ by $X_{\epsilon} = e^{h_{\epsilon}}L_{n} +$ complex tangential terms, and then suitably extended (see p.252). $\rho_{\epsilon}$ is not an actual defining function because it is not real valued (only approximately so). It satisfies
\begin{equation}\label{fix2}
\left |\sum_{j,k = 1}^{n} \frac{\partial^{2}\widehat{\rho_{\epsilon}}}{\partial z_{j}\partial\overline{z_{k}}}(P)w_{j}\overline{w_{k}} \right | = O(\epsilon)|w|^{2},\;
w \in span_{\mathbb{C}}\{N(P), L_{n}(P)\} \; .
\end{equation}
It now suffices to take the family $\rho_{\epsilon} :=$ real part of $(\widehat{\rho_{\epsilon}})$; \eqref{fix2} carries over by taking real and imaginary parts.

\section{On some operators arising from elliptic regularization}

In this section, we give some properties (mainly estimates in Sobolev norms) of operators arising from the regularized $\overline{\partial}$-Neumann problem (\cite{FollandKohn72}, section 3, chapter 2, \cite{Taylor96}, section 5, chapter 12). For $\delta > 0$, $\Box_{\delta,q}$ is the selfadjoint operator defined by the quadratic form
\begin{equation}\label{reg1}
Q_{\delta,q}(u,\overline{u}) = \|\overline{\partial}u\|^{2}
+ \|\overline{\partial}^{*}u\|^{2} + \delta\|\nabla u\|^{2} \;,
\end{equation}
where $\nabla u$ denotes the vector of all (first) derivatives of all components of $u$. The form domain is $W_{(0,q)}^{1}(\Omega) \cap dom(\overline{\partial}^{*})$. $\Box_{\delta,q}$ has a bounded inverse $N_{\delta,q}$. In fact, because $Q_{\delta,q}$ dominates $\|u\|_{1}^{2}$ (for $\delta > 0$), the form is coercive, and the elliptic theory applies (see e,g. \cite{Taylor96}): $N_{\delta,q}$ maps $C^{\infty}_{(0,q)}(\overline{\Omega})$ continuously into itself. Computing $\Box_{\delta,q}$ and the free boundary condition gives

\begin{lemma}\label{Reg1}
Let $u \in dom(\Box_{\delta,q})$, $u = \sum^{\prime}_{J}u_{J}d\overline{z_{J}}$. Then
\begin{equation}\label{reg2}
 \Box_{\delta,q}u = -\sum^{\prime}_{J}(1/4 + \delta)\Delta u_{J}d\overline{z_{J}}\;.
\end{equation}
and
\begin{equation}\label{reg3}
 (\overline{\partial}u)_{norm}(z) + \delta((\partial / \partial\nu)u)_{tan}(z)  = 0 \;, \; z \in  b\Omega \; .
\end{equation}
\end{lemma}
Here $(\overline{\partial}u)_{norm}$ denotes the normal \emph{component} of $\overline{\partial}u$ (a $(0,q)$-form), $((\partial /  \partial\nu)u)_{tan}$ denotes he tangential \emph{part} of $(\partial /\partial\nu)u$ ( also a $q$-form); $(\partial /\partial\nu)u = \sum^{\prime}_{J} (\partial /\partial\nu)u_{J}\;d\overline{z_{J}}$. The lemma is obtained in the same way as the corresponding statements for $\Box_{q}$ (but compare \cite{Taylor96}, p.410). $\Box$

Denote by $P_{q}$ the Bergman projection on $(0,q)$-forms, that is , the orthogonal projection from $\mathcal{L}^{2}_{(0,q)}(\Omega)$ onto the closed subspace of $\overline{\partial}$-closed forms. For $t>0$, denote by $N_{t,q}$ the $\overline{\partial}$-Neumann operator resulting when the $\overline{\partial}$-Neumann problem is set up with respect to the weight factor $w_{t}(z) = e^{-t|z|^{2}}$ (\cite{Kohn73}).

\begin{lemma}\label{Reg2}
Let $u \in \mathcal{L}^{2}_{(0,q)}(\Omega)$, $t>0$, $1 \leq q \leq n$. Then
\begin{equation}\label{reg4}
N_{\delta, q}u = P_{q}w_{t}N_{t,q}\overline{\partial}(w_{-t}\overline{\partial}^{*}N_{\delta,q}u)
+ (Id - P_{q})\overline{\partial}_{t}^{*}N_{t,q+1}\overline{\partial}N_{\delta,q}u \; .
\end{equation}
In particular, if $P_{q}$ satisfies Sobolev estimates in $W^{s}_{(0,q)}(\Omega)$ for some $s>0$, then
\begin{equation}\label{reg5}
\|N_{\delta,q}u\|_{s} \leq C(\|\overline{\partial}N_{\delta,q}u\|_{s} +
\|\overline{\partial}^{*}N_{\delta,q}u\|_{s}) \; .
\end{equation}
\end{lemma}

\emph{Proof}: The proof of the lemma results from the ideas in \cite{BoasStraube90}. We have
\begin{equation}\label{reg6}
N_{\delta,q}u = N_{\delta,q}(\overline{\partial}\overline{\partial}^{*}N_{q}u +
\overline{\partial}^{*}\overline{\partial}N_{q}u) =
(N_{\delta,q}\overline{\partial})(\overline{\partial}^{*}N_{q})u +
(N_{\delta,q}\overline{\partial}^{*})(\overline{\partial}N_{q})u \; .
\end{equation}
Since $N_{\delta,q} = (N_{\delta,q})^{*}$, taking adjoints gives 
\begin{equation}\label{reg7}
N_{\delta,q}u = (\overline{\partial}^{*}N_{q})^{*}(\overline{\partial}^{*}N_{\delta,q})u +
(\overline{\partial}N_{q})^{*}(\overline{\partial}N_{\delta,q})u \; .
\end{equation}
Expressing $\overline{\partial}^{*}N_{q}$ and $\overline{\partial}N_{q} = N_{q+1}\overline{\partial}$ in terms of weighted operators, as in \cite{BoasStraube90}, gives
\begin{equation}\label{reg8}
N_{\delta,q}u = P_{q}w_{t}N_{t,q}\overline{\partial}\left (w_{-t}(Id - P_{q-1})\overline{\partial}^{*}N_{\delta,q}u\right ) +
(Id - P_{q})\overline{\partial}_{t}^{*}N_{t,q+1}P_{q+1}\overline{\partial}N_{\delta,q}u \; .
\end{equation}
This is \eqref{reg4}, because $(Id - P_{q-1})\overline{\partial}^{*}N_{\delta,q}u =
\overline{\partial}^{*}N_{\delta,q}u$ and $P_{q+1}\overline{\partial}N_{\delta,q}u = 
\overline{\partial}N_{\delta,q}u$. \eqref{reg5} is a consequence of \eqref{reg4} and Kohn's weighted theory (\cite{Kohn73}): for a given $s \geq 0$, we may take $t$ big enough so that both $N_{t,q}\overline{\partial}$ and $\overline{\partial}_{t}^{*}N_{t,q+1}$ are continuous in $W^{s}$ (see \cite{BoasStraube91}, p.84-85 for details). $\Box$

\begin{lemma}\label{Reg3}
Let $k \in \mathbb{N}$. There is a constant $C = C(k)$ such that when $\delta >0$ and $u \in C^{\infty}_{(0,q)}(\overline{\Omega})$
\begin{equation}\label{reg9}
\delta^{2}\|N_{\delta,q}u\|_{k+2}^{2} \leq C(\|u\|_{k}^{2} + \delta\|\nabla N_{\delta,q}u\|_{k}^{2} + \delta^{2}\|u\|_{k+1}^{2}) \; .
\end{equation}
\end{lemma}

\emph{Proof}: We use $|||u|||_{k}$ to denote tangential Sobolev norms, and we denote by $\Lambda_{s}$ the standard tangential operators of order $s$, see e.g. \cite{FollandKohn72}, chapter 2, section 4, or \cite{ChenShaw01}, section 5.2. We use that $\|N_{\delta,q}u\|_{k+2}^{2} \leq C(|||N_{\delta,q}u|||_{k+2}^{2} + \|u\|_{k+1}^{2})$; see for example \cite{ChenShaw01}, Lemma 5.2.4. The lemma is stated for $\Box_{q}$ (i.e.$N_{q}$), but as the authors point out (p.102), one can repeat the proof for $\Box_{\delta,q}$. Now
\begin{multline}\label{reg10}
\delta^{2}|||N_{\delta,q}u|||_{k+2}^{2} \leq C\delta^{2}\|\Lambda_{k+1}N_{\delta,q}u\|_{1}^{2}
\leq C\delta Q_{\delta}(\Lambda^{k+1}N_{\delta,q}u, \Lambda^{k+1}N_{\delta,q}u) \\
\leq C\delta\left ( |Q_{\delta}(N_{\delta,q}u, (\Lambda^{k+1})^{*}\Lambda^{k+1}N_{\delta,q}u) |
+ O(|||\nabla N_{\delta,q}u|||_{k}^{2}) \right ) \; .
\end{multline}
The last inequality in \eqref{reg10} comes from Lemma 3.1 in \cite{KohnNirenberg65}, see also Lemma 2.4.2 in \cite{FollandKohn72}. The first term on the right hand side in \eqref{reg10} equals, after moving $k$ factors $\Lambda^{*}$ back to the left
\begin{multline}\label{reg11}
C\delta\left ( |(\Lambda^{k}u, \Lambda^{*}\Lambda^{k+1}N_{\delta,q}u) | +
O(|||\nabla N_{\delta,q}u|||_{k}^{2})\right ) \\
\leq C\delta\|u\|_{k}|||N_{\delta,q}u|||_{k+2}
+ \delta O(|||\nabla N_{\delta,q}u|||_{k}^{2}) \; .
\end{multline}
We estimate the first term on the right hand side in \eqref{reg11} as $(s.c.)\delta^{2}|||N_{\delta,q}u|||_{k+2}^{2} + (l.c.)\|u\|_{k}^{2}$; absorbing $(s.c.)\delta^{2}|||N_{\delta,q}u|||_{k+2}^{2}$ completes the proof of Lemma \ref{Reg3}. $\Box$

\begin{lemma}\label{Reg4}
Let $k \in \mathbb{N}$. Then we have the estimate
\begin{multline}\label{reg12}
\|\overline{\partial}\overline{\partial}^{*}N_{\delta, q}u\|_{k}^{2} +
\|\vartheta\overline{\partial}N_{\delta, q}u\|_{k}^{2} \\
\leq C\left (\|\overline{\partial}N_{\delta, q}u\|_{k}^{2} + \|u\|_{k}^{2} + \|u\|_{1}^{2} +
\delta\|\nabla N_{\delta, q}u\|_{k}^{2} + \delta^{2}\|u\|_{k+1}^{2} \right ) \; ,
\end{multline}
with a constant $C$ independent of $\delta$ (and of course of $u$).
\end{lemma}

We remark that the term $\|u\|_{1}^{2}$ is only relevant when $k=0$. It arises in connection with trace theorems for functions which are only in $\mathcal{L}^{2}(\Omega)$ (see below). Since $\overline{\partial}\overline{\partial}^{*}N_{\delta, q}u + \vartheta\overline{\partial}N_{\delta, q}u =\left (1/(1+4\delta)\right )\Box_{\delta,q}N_{\delta,q}u =
\left (1/(1+4\delta)\right )u$, it suffices to estimate one of the two terms in \eqref{reg12}, say $\|\vartheta\overline{\partial}N_{\delta, q}u\|_{k}^{2}$. We have
\begin{equation}\label{reg13}
\|\vartheta\overline{\partial}N_{\delta, q}u\|_{k}^{2} 
\leq C\left (|||\vartheta\overline{\partial}N_{\delta, q}u|||_{k}^{2} +
\|\overline{\partial}\vartheta\overline{\partial}N_{\delta, q}u\|_{k-1}^{2} +
\|\vartheta\overline{\partial}N_{\delta, q}u\|_{k-1}^{2} \right ) \; .
\end{equation}
Note that $\overline{\partial}\vartheta\overline{\partial}N_{\delta, q}u =
\overline{\partial}(\vartheta\overline{\partial}+\overline{\partial}\overline{\partial}^{*})N_{\delta, q}u = \left (1/(1+4\delta)\right )\overline{\partial}u$; therefore, the middle term on the right hand side of \eqref{reg13} is of order $\|u\|_{k}^{2}$. For the first term on the right hand side of \eqref{reg13} we have, denoting by $T^{k}$ a tangential differential operator of order $k$,
\begin{multline}\label{reg14}
(T^{k}\vartheta\overline{\partial}N_{\delta, q}u, T^{k}\vartheta\overline{\partial}N_{\delta, q}u) \\
= (\vartheta\overline{\partial}N_{\delta, q}u, T^{2k}\vartheta\overline{\partial}N_{\delta, q}u) + 
O(\|\vartheta\overline{\partial}N_{\delta, q}u\|_{k-1}\|\vartheta\overline{\partial}N_{\delta, q}u\|_{k} ) \; .
\end{multline}
The main term in \eqref{reg14} equals
\begin{multline}\label{reg15}
(\overline{\partial}N_{\delta,q}u, \overline{\partial}T^{2k}\vartheta\overline{\partial}N_{\delta, q}u ) -
\int_{b\Omega}\langle(\overline{\partial}N_{\delta,q}u)_{n}, T^{2k}\vartheta\overline{\partial}N_{\delta, q}u \rangle \\
= (\overline{\partial}N_{\delta,q}u, T^{2k}\overline{\partial}\vartheta\overline{\partial}N_{\delta, q}u )
+ (\overline{\partial}N_{\delta,q}u, [\overline{\partial}, T^{2k}]\vartheta\overline{\partial}N_{\delta, q}u ) \\
+ \int_{b\Omega}\langle\delta ((\partial / \partial\nu)N_{\delta,q}u)_{tan}, T^{2k}\vartheta\overline{\partial}N_{\delta, q}u \rangle \; .
\end{multline}
In the second term on the right hand side of \eqref{reg15}, we expand the commutator as in section 4, to get a main term of order $\|\overline{\partial}N_{\delta,q}u\|_{k}\|\vartheta\overline{\partial}N_{\delta,q}u\|_{k}$. The first term on the right hand side of \eqref{reg15} equals (up to a constant) $(\overline{\partial}N_{\delta,q}u, T^{2k}\overline{\partial}u)$. It is estimated by performing the above computations in reverse order to arrive at $\|\vartheta\overline{\partial}N_{\delta,q}u\|_{k}\|u\|_{k}$ (for the main term). To estimate the boundary integral, we use duality of Sobolev spaces on the boundary and the trace theorem in $W_{(0,q)}^{k}(\Omega)$ and $W_{(0,q)}^{k+1}(\Omega)$, respectively, to bound this term by$\|\delta (\partial /\partial\nu)N_{\delta,q}u\|_{k+1}\|\vartheta\overline{\partial}N_{\delta,q}u\|_{k}$. It is here that we need $k \geq 1$. Using Lemma \ref{Reg3}, we can estimate 
\begin{equation}\label{reg16}
\|\delta (\partial / \partial\nu)N_{\delta,q}u\|_{k+1}^{2} \leq \delta^{2}\|N_{\delta,q}u\|_{k+2}^{2} \leq C(\|u\|_{k}^{2} + \delta\|\nabla N_{\delta,q}u\|_{k}^{2} + \delta^{2}\|u\|_{k+1}^{2} ) \; .
\end{equation}
The third term on the right hand side of \eqref{reg13} is dominated by $(s.c.)\|\vartheta\overline{\partial}N_{\delta,q}u\|_{k}^{2} + (l.c.)\|\vartheta\overline{\partial}N_{\delta,q}u\|_{-1}^{2}$, by interpolation of Sobolev norms. The first term can be absorbed, the second is dominated by $\|\overline{\partial}N_{\delta,q}u\|_{0}^{2} \leq \|u\|^{2}$.

When $k=0$, the above integration by parts argument still works, but the trace estimate from $\mathcal{L}^{2}(\Omega)$ to $W^{-1/2}(\Omega)$ applied to $\vartheta\overline{\partial}N_{\delta,q}u$ now needs the term $\|\Delta\vartheta\overline{\partial}N_{\delta,q}u\|_{-1} = \left (4/ (1+4\delta)\right )\|\vartheta\overline{\partial}u\|_{-1} \leq C\|u\|_{1}$ (see e.g. \cite{LionsMagenes72}). (We are not striving for an optimal estimate for this case, but simply one that suffices for our purposes.) 

This completes the proof of Lemma \ref{Reg4}. $\Box$


%
%

\providecommand{\bysame}{\leavevmode\hbox to3em{\hrulefill}\thinspace}

\end{document}